\documentclass[10pt]{article}
\usepackage{amsmath}
\usepackage{latexsym}
\usepackage{amssymb}
\usepackage{amsfonts}
\usepackage{amsthm}
\title{MINIMUM PARAMETRIZATION OF \\
THE CAUCHY STRESS OPERATOR}
\author{J.-F. POMMARET \\ CERMICS, Ecole des Ponts ParisTech\\ 
 jean-francois.pommaret@wanadoo.fr \\
ORCID: 0000-0003-0907-2601 }
\date{  }
\begin{document}
\maketitle

\noindent
{\bf ABSTRACT}    \\

\noindent
When ${\cal{D}}:\xi \rightarrow \eta$ is a linear differential operator, a "direct problem " is to find the generating compatibility conditions (CC) in the form of an operator ${\cal{D}}_1:\eta \rightarrow \zeta$ such that ${\cal{D}}\xi=\eta$ implies ${\cal{D}}_1\eta=0$. When ${\cal{D}}$ is involutive, the procedure provides successive first order involutive operators ${\cal{D}}_1, ... , {\cal{D}}_n$ when the ground manifold has dimension $n$. Conversely, when ${\cal{D}}_1$ is given, a more difficult " inverse problem " is to look for an operator ${\cal{D}}: \xi \rightarrow \eta$  having the generating CC ${\cal{D}}_1\eta=0$. If this is possible, that is when the differential module defined by ${\cal{D}}_1$ is torsion-free, one shall say that the operator ${\cal{D}}_1$ is parametrized by ${\cal{D}}$ and there is no relation in general between ${\cal{D}}$ and ${\cal{D}}_2$. The parametrization is said to be " minimum "  if the differential module defined by ${\cal{D}}$ has a vanishing differential rank and is thus a torsion module. The parametrization of the Cauchy stress operator in arbitrary dimension $n$ has attracted many famous scientists (G.B. Airy in 1863 for $n=2$, J.C. Maxwell in 1863, G. Morera and E. Beltrami in 1892 for $n=3$, A. Einstein in 1915 for $n=4$) . The aim of this paper is to prove that all these works are explicitly using the Einstein operator (which cannot be parametrized) and not the Ricci operator. As a byproduct, they are all based on a confusion between the so-called $div$ operator induced from the Bianchi operator ${\cal{D}}_2$ and the Cauchy operator which is the formal adjoint of the Killing operator ${\cal{D}}$ parametrizing the Riemann operator ${\cal{D}}_1$ for an arbitrary $n$. We also present the similar situation met in the study of contact structures when $n=3$. Like the Michelson and Morley experiment, it is an open historical problem to know whether Einstein was aware of these  previous works or not, as the comparison needs no comment.

 \vspace{4cm}


\newpage

\noindent
{\bf 1) INTRODUCTION}  \\

We start recalling the basic tools from the formal theory of systems of partial differential (PD) equations and differential modules needed in order to understand and solve the parametrization problem presented in the abstract. Then we provide the example of the system of infinitesimal Lie equations defining contact transformations and conclude the paper with the general parametrization problem existing in continuum mechanics for an arbitrary dimension of the ground manifold. As these new tools are difficult and not so well known, we advise the interested reader to follow them step by step on the explicit motivating examples illustrating this paper.  \\

\noindent
{\bf A) SYSTEM THEORY}:  \\

If $X$ is a manifold of dimension $n$ with local coordinates $(x)=(x^1, ... ,x^n)$, we denote as usual by $T=T(X)$ the {\it tangent bundle} of $X$, by $T^*=T^*(X)$ the {\it cotangent bundle}, by ${\wedge}^rT^*$ the {\it bundle of r-forms} and by $S_qT^*$ the {\it bundle of q-symmetric tensors}. More generally, let $E$ be a {\it vector bundle} over $X$ with local coordinates $(x^i,y^k)$ for $i=1,...,n$ and $k=1,...,m$ simply denoted by $(x,y)$, {\it projection} $\pi:E\rightarrow X:(x,y)\rightarrow (x)$ and changes of local coordinate $\bar{x}=\varphi (x), \bar{y}=A(x)y$. We shall denote by $E^*$ the vector bundle obtained by inverting the matrix $A$ of the changes of coordinates , exactly like $T^*$ is obtained from $T$. We denote by $f:X\rightarrow E: (x)\rightarrow (x,y=f(x))$ a global {\it section} of $E$, that is a map such that $\pi\circ f=id_X$ but local sections over an open set $U\subset X$ may also be considered when needed. Under a change of coordinates, a section transforms like $\bar{f}(\varphi(x))=A(x)f(x)$ and the changes of the derivatives can also be obtained with more work. We shall denote by $J_q(E)$ the {\it q-jet bundle} of $E$ with local coordinates $(x^i, y^k, y^k_i, y^k_{ij},...)=(x,y_q)$ called {\it jet coordinates} and sections $f_q:(x)\rightarrow (x,f^k(x), f^k_i(x), f^k_{ij}(x), ...)=(x,f_q(x))$ transforming like the sections $j_q(f):(x) \rightarrow (x,f^k(x), {\partial}_if^k(x), {\partial}_{ij}f^k(x), ...)=(x,j_q(f)(x))$ where both $f_q$ and $j_q(f)$ are over the section $f$ of $E$. For any $q\geq 0$, $J_q(E)$ is a vector bundle over $X$ with projection ${\pi}_q$ while $J_{q+r}(E)$ is a vector bundle over $J_q(E)$ with projection ${\pi}^{q+r}_q, \forall r\geq 0$.\\

\noindent
{\bf DEFINITION  1.A.1}: A {\it linear system} of order $q$ on $E$ is a vector sub-bundle $R_q\subset J_q(E)$ and a {\it solution} of $R_q$ is a section $f$ of $E$ such that $j_q(f)$ is a section of $R_q$. With a slight abuse of language, the set of local solutions will be denoted by $\Theta \subset E$.\\
  
Let $\mu=({\mu}_1,...,{\mu}_n)$ be a multi-index with {\it length} ${\mid}\mu{\mid}={\mu}_1+...+{\mu}_n$, {\it class} $i$ if ${\mu}_1=...={\mu}_{i-1}=0,{\mu}_i\neq 0$ and $\mu +1_i=({\mu}_1,...,{\mu}_{i-1},{\mu}_i +1, {\mu}_{i+1},...,{\mu}_n)$. We set $y_q=\{y^k_{\mu}{\mid} 1\leq k\leq m, 0\leq {\mid}\mu{\mid}\leq q\}$ with $y^k_{\mu}=y^k$ when ${\mid}\mu{\mid}=0$. If $E$ is a vector bundle over $X$ and $J_q(E)$ is the $q$-{\it jet bundle} of $E$, then both sections $f_q\in J_q(E)$ and $j_q(f)\in J_q(E)$ are over the section $f\in E$. There is a natural way to distinguish them by introducing the {\it Spencer}  operator $d:J_{q+1}(E)\rightarrow T^*\otimes J_q(E)$ with components $(df_{q+1})^k_{\mu,i}(x)={\partial}_if^k_{\mu}(x)-f^k_{\mu+1_i}(x)$. The kernel of $d$ consists of sections such that $f_{q+1}=j_1(f_q)=j_2(f_{q-1})=...=j_{q+1}(f)$. Finally, if $R_q\subset J_q(E)$ is a {\it system} of order $q$ on $E$ locally defined by linear equations ${\Phi}^{\tau}(x,y_q)\equiv a^{\tau\mu}_k(x)y^k_{\mu}=0$ and local coordinates $(x,z)$ for the parametric jets up to order $q$, the $r$-{\it prolongation} $R_{q+r}={\rho}_r(R_q)=J_r(R_q)\cap J_{q+r}(E)\subset J_r(J_q(E))$ is locally defined when $r=1$ by the linear equations ${\Phi}^{\tau}(x,y_q)=0, d_i{\Phi}^{\tau}(x,y_{q+1})\equiv a^{\tau\mu}_k(x)y^k_{\mu+1_i}+{\partial}_ia^{\tau\mu}_k(x)y^k_{\mu}=0$ and has {\it symbol} $g_{q+r}=R_{q+r}\cap S_{q+r}T^*\otimes E\subset J_{q+r}(E)$ if one looks at the {\it top order terms}. If $f_{q+1}\in R_{q+1}$ is over $f_q\in R_q$, differentiating the identity $a^{\tau\mu}_k(x)f^k_{\mu}(x)\equiv 0$ with respect to $x^i$ and substracting the identity $a^{\tau\mu}_k(x)f^k_{\mu+1_i}(x)+{\partial}_ia^{\tau\mu}_k(x)f^k_{\mu}(x)\equiv 0$, we obtain the identity $a^{\tau\mu}_k(x)({\partial}_if^k_{\mu}(x)-f^k_{\mu+1_i}(x))\equiv 0$ and thus the restriction $d:R_{q+1}\rightarrow T^*\otimes R_q$. More generally, we have the restriction:   \\
\[   d: {\wedge}^sT^* \otimes R_{q+1} \rightarrow {\wedge}^{s+1}T^* \otimes R_q: (f^k_{\mu,i}(x)dx^I) \rightarrow (({\partial}_if^k_{\mu,i}(x) - f^k_{\mu + 1_i,I}(x))dx^i \wedge dx^I)\] 
with standard multi-index notation for exterior forms and one can easily check that $d\circ d=0$.\\
The restriction of $-d$ to the symbol is called the {\it Spencer} map $\delta: {\wedge}^sT^*\otimes g_{q+1} \rightarrow {\wedge}^{s+1}T^* \otimes g_q $ and $\delta \circ \delta=0$ similarly ([22-25],[28],[41],[49]). \\
    
\noindent
{\bf DEFINITION 1.A.2}: A system $R_q$ is said to be {\it formally integrable} when all the equations of order $q+r$ are obtained by $r$ prolongations {\it only}, $\forall r\geq 0$ or, equivalently, when the projections ${\pi}^{q+r+s}_{q+r}:R_{q+r+s}\rightarrow R{q+r}$ are epimorphisms $\forall r,s\geq0$.  \\

Finding an intrinsic test has been achieved by D.C. Spencer in 1970 ([49]) along coordinate dependent lines sketched by M. Janet in 1920 ([11]). The next procedure providing a {\it Pommaret basis} and where  {\it one may have to change linearly the independent variables if necessary}, is intrinsic even though it must be checked in a particular coordinate system called $\delta$-{\it regular} ([22],[25],[29]).  \\

$\bullet$ {\it Equations of class} $n$: Solve the maximum number ${\beta}^n_q$ of equations with respect to the jets of order $q$ and class $n$. Then call $(x^1,...,x^n)$ {\it multiplicative variables}.\\

$\bullet$ {\it Equations of class} $i\geq 1$: Solve the maximum number ${\beta}^i_q$ of {\it remaining} equations with respect to the jets of order $q$ and class $i$. Then call $(x^1,...,x^i)$ {\it multiplicative variables} and $(x^{i+1},...,x^n)$ {\it non-multiplicative variables}.\\

$\bullet$ {\it Remaining equations equations of order} $\leq q-1$: Call $(x^1,...,x^n)$ {\it non-multiplicative variables}.\\

\noindent
In actual practice, we shall use a {\it Janet tabular} where the multiplicative "variables" are in upper left position while the non-multiplicative variables are represented by dots in lower right position.  \\

\noindent
{\bf DEFINITION 1.A 3}: A system of PD equations is said to be {\it involutive} if its first prolongation can be obtained by prolonging its equations only with respect to the corresponding multiplicative variables. In that case, we may introduce the {\it characters} ${\alpha}^i_q=m\frac{(q+n-i-1)!}{(q-1)!((n-i)!}-{\beta}^i_q$ for $i=1, ..., n$ with ${\alpha}^1_q\geq ... \geq {\alpha}^n_q\geq 0$ and we have $dim(g_q)={\alpha}^1_q+...+ {\alpha}^n_q$ while $dim(g_{q+1})={\alpha}^1_q+...+ n {\alpha}^n_q$. \\

\noindent
{\bf REMARK 1.A.4}: As long as the {\it Prolongation}/{\it Projection} (PP) procedure has not been achieved in order to get an involutive system, {\it nothing} can be said about the CC (Fine examples can be found in [41] and the recent [45]).  \\

\noindent
{\bf REMARK  1.A.5}: A proof that the second order system defined by Einstein equations is involutive has been given by J. Gasqui in $1982$ but this paper cannot be applied to the minimum parametrizations that need specific $\delta$-regular coordinates as we shall see ([8]).  \\

When $R_q$ is involutive, the linear differential operator ${\cal{D}}:E\stackrel{j_q}{\rightarrow} J_q(E)\stackrel{\Phi}{\rightarrow} J_q(E)/R_q=F_0$ of order $q$ is said to be {\it involutive} and its space of solutions is defined by the kernel exact sequence $0 \rightarrow \Theta \rightarrow E \stackrel{{\cal{}}}{\longrightarrow} F_0$. One has the canonical {\it linear Janet sequence} (Introduced in [19]):\\
\[  0 \longrightarrow  \Theta \longrightarrow E \stackrel{\cal{D}}{\longrightarrow} F_0 \stackrel{{\cal{D}}_1}{\longrightarrow}F_1 \stackrel{{\cal{D}}_2}{\longrightarrow} ... \stackrel{{\cal{D}}_n}{\longrightarrow} F_n \longrightarrow 0   \]
where each other operator is first order involutive and generates the {\it compatibility conditions} (CC) of the preceding one. Similarly, introducing the Spencer bundles $C_r= {\wedge}^rT^* \otimes R_q / \delta ({\wedge}^{r-1}T^* \otimes g_{q+1})$ we obtain the canonical {\it linear Spencer sequence} induced by the Spencer operator:  \\
\[   0\longrightarrow \Theta \stackrel{j_q}{\longrightarrow} C_0 \stackrel{D_1}{\longrightarrow} C_1 \stackrel{D_2}{\longrightarrow}... \stackrel{D_n}{\longrightarrow} C_n \longrightarrow  0  \]

\noindent
{\bf B) MODULE THEORY}:   \\

Let $K$ be a {\it differential field} with $n$ commuting derivations $({\partial}_1,...,{\partial}_n)$ and consider the ring $D=K[d_1,...,d_n]=K[d]$ of differential operators with coefficients in $K$ with $n$ commuting formal derivatives satisfying $d_ia=ad_i + {\partial}_ia$ in the operator sense. If $P=a^{\mu}d_{\mu}\in D=K[d]$, the highest value of ${\mid}\mu {\mid}$ with $a^{\mu}\neq 0$ is called the {\it order} of the {\it operator} $P$ and the ring $D$ with multiplication $(P,Q)\longrightarrow P\circ Q=PQ$ is filtred by the order $q$ of the operators. We have the {\it filtration} $0\subset K=D_0\subset D_1\subset  ... \subset D_q \subset ... \subset D_{\infty}=D$. As an algebra, $D$ is generated by $K=D_0$ and $T=D_1/D_0$ with $D_1=K\oplus T$ if we identify an element $\xi={\xi}^id_i\in T$ with the vector field $\xi={\xi}^i(x){\partial}_i$ of differential geometry, but with ${\xi}^i\in K$ now. It follows that $D={ }_DD_D$ is a {\it bimodule} over itself, being at the same time a left $D$-module by the composition $P \longrightarrow QP$ and a right $D$-module by the composition $P \longrightarrow PQ$. We define the {\it adjoint} functor $ad:D \longrightarrow D^{op}:P=a^{\mu}d_{\mu} \longrightarrow  ad(P)=(-1)^{\mid \mu \mid}d_{\mu}a^{\mu}$ and we have $ad(ad(P))=P$ both with $ad(PQ)=ad(Q)ad(P), \forall P,Q\in D$. Such a definition can be extended to any matrix of operators by using the transposed matrix of adjoint operators (See [4],[12],[28],[29],[34],[38],[48] for more details and applications to control theory or mathematical physics). \\

Accordingly, if $y=(y^1, ... ,y^m)$ are differential indeterminates, then $D$ acts on $y^k$ by setting $d_iy^k=y^k_i \longrightarrow d_{\mu}y^k=y^k_{\mu}$ with $d_iy^k_{\mu}=y^k_{\mu+1_i}$ and $y^k_0=y^k$. We may therefore use the jet coordinates in a formal way as in the previous section. Therefore, if a system of OD/PD equations is written in the form ${\Phi}^{\tau}\equiv a^{\tau\mu}_ky^k_{\mu}=0$ with coefficients $a\in K$, we may introduce the free differential module $Dy=Dy^1+ ... +Dy^m\simeq D^m$ and consider the differential {\it module of equations} $I=D\Phi\subset Dy$, both with the residual {\it differential module} $M=Dy/D\Phi$ or $D$-module and we may set $M={ }_DM$ if we want to specify the ring of differential operators. We may introduce the formal {\it prolongation} with respect to $d_i$ by setting $d_i{\Phi}^{\tau}\equiv a^{\tau\mu}_ky^k_{\mu+1_i}+({\partial}_ia^{\tau\mu}_k)y^k_{\mu}$ in order to induce maps $d_i:M \longrightarrow M:{\bar{y} }^k_{\mu} \longrightarrow {\bar{y}}^k_{\mu+1_i}$ by residue with respect to $I$ if we use to denote the residue $Dy \longrightarrow M: y^k \longrightarrow {\bar{y}}^k$ by a bar like in algebraic geometry. However, for simplicity, we shall not write down the bar when the background will indicate clearly if we are in $Dy$ or in $M$. As a byproduct, the differential modules we shall consider will always be {\it finitely generated} ($k=1,...,m<\infty$) and {\it finitely presented} ($\tau=1, ... ,p<\infty$). Equivalently, introducing the {\it matrix of operators} ${\cal{D}}=(a^{\tau\mu}_kd_{\mu})$ with $m$ columns and $p$ rows, we may introduce the morphism $D^p \stackrel{{\cal{D}}}{\longrightarrow} D^m:(P_{\tau}) \longrightarrow (P_{\tau}{\Phi}^{\tau})$ over $D$ by acting with $D$ {\it on the left of these row vectors} while acting with ${\cal{D}}$ {\it on the right of these row vectors} by composition of operators with $im({\cal{D}})=I$. The {\it presentation} of $M$ is defined by the exact cokernel sequence $D^p \stackrel{{\cal{D}}}{\longrightarrow} D^m \longrightarrow M \longrightarrow 0 $. We notice that the presentation only depends on $K, D$ and $\Phi$ or $ \cal{D}$, that is to say never refers to the concept of (explicit local or formal) solutions. It follows from its definition that $M$ can be endowed with a {\it quotient filtration} obtained from that of $D^m$ which is defined by the order of the jet coordinates $y_q$ in $D_qy$. We have therefore the {\it inductive limit} $0 \subseteq M_0 \subseteq M_1 \subseteq ... \subseteq M_q \subseteq ... \subseteq M_{\infty}=M$ with $d_iM_q\subseteq M_{q+1}$ and $M=DM_q$ for $q\gg 0$ with prolongations $D_rM_q\subseteq M_{q+r}, \forall q,r\geq 0$.  \\

\noindent
{\bf DEFINITION 1.B.1}: An exact sequence of morphisms finishing at $M$ is said to be a {\it resolution} of $M$. If the differential modules involved apart from $M$ are free, that is isomorphic to a certain power of $D$, we shall say that we have a {\it free resolution} of $M$.  \\

Having in mind that $K$ is a left $D$-module with the action $(D,K) \longrightarrow K:(d_i,a)\longrightarrow {\partial}_ia$ and that $D$ is a bimodule over itself, {\it we have only two possible constructions}:  \\

\noindent
{\bf DEFINITION 1.B.2}: We may define the right differential module $hom_D(M,D)$.  \\

\noindent
{\bf DEFINITION 1.B.3}: We define the {\it system} $R=hom_K(M,K)$ and set $R_q=hom_K(M_q,K)$ as the {\it system of order} $q$. We have the {\it projective limit} $R=R_{\infty} \longrightarrow ... \longrightarrow R_q \longrightarrow ... \longrightarrow R_1 \longrightarrow R_0$. It follows that $f_q\in R_q:y^k_{\mu} \longrightarrow f^k_{\mu}\in K$ with $a^{\tau\mu}_kf^k_{\mu}=0$ defines a {\it section at order} $q$ and we may set $f_{\infty}=f\in R$ for a {\it section} of $R$. For an arbitrary differential field $K$, {\it such a definition has nothing to do with the concept of a formal power series solution} ({\it care}).\\

\noindent
{\bf PROPOSITION 1.B.4}: When $M$ is a left $D$-module, then $R$ is also a left $D$-module. \\

\noindent
{\it Proof}: As $D$ is generated by $K$ and $T$ as we already said, let us define:  \\
\[  (af)(m)=af(m), \hspace{4mm} \forall a\in K, \forall m\in M \]
\[ (\xi f)(m)=\xi f(m)-f(\xi m), \hspace{4mm} \forall \xi=a^id_i\in T,\forall m\in M  \]
In the operator sense, it is easy to check that $d_ia=ad_i+{\partial}_ia$ and that $\xi\eta - \eta\xi=[\xi,\eta]$ is the standard bracket of vector fields. We finally 
get $(d_if)^k_{\mu}=(d_if)(y^k_{\mu})={\partial}_if^k_{\mu}-f^k_{\mu +1_i}$ and thus recover {\it exactly} the Spencer operator of the previous section though {\it this is not evident at all}. We also get $(d_id_jf)^k_{\mu}={\partial}_{ij}f^k_{\mu}-{\partial}_if^k_{\mu+1_j}-{\partial}_jf^k_{\mu+1_i}+f^k_{\mu+1_i+1_j} \Longrightarrow d_id_j=d_jd_i, \forall i,j=1,...,n$ and thus $d_iR_{q+1}\subseteq R_q\Longrightarrow d_iR\subset R$ induces a well defined operator $R\longrightarrow T^*\otimes R:f \longrightarrow dx^i\otimes d_if$. This operator has been first introduced, up to sign, by F.S. Macaulay as early as in $1916$ but this is still not ackowledged ([17]). For more details on the Spencer operator and its applications, the reader may look at ([31],[37-39],[42]).  \\
\hspace*{12cm}   Q.E.D.  \\

\noindent
{\bf DEFINITION 1.B.5}: With any differential module $M$ we shall associate the {\it graded module} $G=gr(M)$ over the polynomial ring $gr(D)\simeq K[\chi]$ by setting $G={\oplus}^{\infty}_{q=0} G_q$ with $G_q=M_q/M_{q+1}$ and we get $g_q=G_q^*$ where the {\it symbol} $g_q$ is defined by the short exact sequences: \\
\[ 0\longrightarrow M_{q-1}\longrightarrow M_q \longrightarrow G_q \longrightarrow 0  \hspace{4mm}  \Longrightarrow \hspace{4mm}  0 \longrightarrow g_q \longrightarrow R_q \longrightarrow R_{q-1} \longrightarrow 0  \]
We have the short exact sequences $0\longrightarrow D_{q-1} \longrightarrow D_q \longrightarrow S_qT \longrightarrow 0 $ leading to $gr_q(D)\simeq S_qT$ and we may set as usual $T^*=hom_K(T,K)$ in a coherent way with differential geometry. \\

The two following definitions, which are well known in commutative algebra, are also valid (with more work) in the case of differential modules (See [28] for more details or the references [9],[21],[29],[47] for an introduction to homological algebra and diagram chasing).  \\

\noindent
{\bf DEFINITION 1.B.6}: The set of elements $t(M) = \{ m \in M \mid \exists 0 \neq P \in D, Pm=0\}\subseteq M$ is a differential module called the {\it torsion submodule} of $M$. More generally, a module $M$ is called a {\it torsion module} if $t(M)=M$ and a {\it torsion-free module} if $t(M)=0$. In the short exact sequence $0 \rightarrow t(M) \rightarrow M \rightarrow M' \rightarrow 0$, the module $M'$ is torsion-free. Its defining module of equations $I'$ is obtained by adding to $I$ a representative basis of $t(M)$ set up to zero and we have thus $I \subseteq I'$.  \\

\noindent
{\bf DEFINITION 1.B.7}: A differential module $F$ is said to be {\it free} if $F \simeq D^r$ for some integer 
$r > 0$ and we shall {\it define} $rk_D(F)=r$. If $F$ is the biggest free dfferential module contained in $M$, then $M/F$ is a torsion differential module and $hom_D(M/F,D)=0$. In that case, we shall {\it define} the {\it differential rank} of $M$ to be $rk_D(M)=rk_D(F)=r$. \\

\noindent
{\bf PROPOSITION 1.B.8}: If $0 \rightarrow M' \rightarrow M \rightarrow M" \rightarrow 0$ is a short exact sequence of differential modules and maps or operators, we have $rk_D(M)=rk_D(M') + rk_D(M")$.  \\

In the general situation, let us consider the sequence $ M' \stackrel{f}{\longrightarrow} M  \stackrel{g}{\longrightarrow} M" $ of modules which may not be exact and define $B =im(f) \subseteq Z = ker(g) \Rightarrow   H=Z/B$. \\

\noindent
{\bf LEMMA  1.B.9}: The kernel of the induced epimorphism $coker(f) \rightarrow coim(g)$ is isomorphic to $H$.  \\

\noindent
{\it Proof}: It follows from a snake chase in the commutative and exact diagram where $coim(g)\simeq im(g)$:  \\
\[  \begin{array}{rcccccccl}
     &  &  &  &  &  &  0  &  &   \\
          &  &  &  &  &  & \downarrow  &  &  \\
     &   &          0     &   &    0   &    &   H   &   &   \\
   &   &  \downarrow  &  &  \downarrow  &   &  \downarrow  &  &  \\

0 & \rightarrow  & B &  \rightarrow & M  & \longrightarrow &   coker(f) & \rightarrow &  0  \\
    &                      &  \downarrow  &    &  \parallel      & &  \downarrow  &   &  \\
0    &  \rightarrow &  Z  &  \rightarrow & M  &  \stackrel{g}{\longrightarrow}  &  coim(g) &  \rightarrow  &  0  \\
        &       &   \downarrow  &  &  \downarrow  &  &  \downarrow   &   &  \\
        &    &       H   &   &   0   &   &  0    &   &   \\
         &  &    \downarrow   &  &  &  &  &  &\\
         &  &   0  &  &  &  &  &  &
\end{array}   \]
\hspace*{12cm}   Q.E.D.   \\

In order to conclude this section, we may say that the main difficulty met when passing from the differential framework to the algebraic framework is the " {\it inversion} " of arrows. Indeed, when an operator is injective, that is when we have the exact sequence $ 0 \rightarrow E \stackrel{{\cal{D}}}{\longrightarrow} F$ with $dim(E)=m, dim(F)=p$, like in the case of the operator $0 \rightarrow E \stackrel{j_q}{\longrightarrow} J_q(E) $, on the contrary, using differenial modules, we have the epimorphism 
$D^p \stackrel{{\cal{D}}}{\longrightarrow} D^m \rightarrow 0$. The case of a formally surjective operator, like the $div$ operator, described by the exact sequence $E \stackrel{{\cal{D}}}{\longrightarrow} F \rightarrow 0$ is now providing the exact sequence of differential modules $ 0 \rightarrow D^p \stackrel{{\cal{D}}}{\longrightarrow} D^m \rightarrow M  \rightarrow 0$ because ${\cal{D}}$ has no CC.  \\

\noindent
{\bf 2) PARAMETRIZATION PROBLEM}     \\

In this section, we shall set up and solve the minimum parametrization problem by comparing the differential geometric approach and the differential algebraic approach. In fact, both sides are essential because certain concepts, like " torsion ", are simpler in the module approach while others, like " involution " are simpler in the opertor approach. However, the reader must never forget that the " {\it extension modules} " or the " {\it side changing functor} " are pure product of {\it differential homological algebra} with no system counterpart. Also, the close link existing between " {\it differential duality} " and " {\it adjoint operator} " may not be evident at all, even for people quite familiar with mathematical physics ([4],[28],[38]).  \\

Let us start with a given linear differential operator $\eta \stackrel{{\cal{D}}_1}{\longrightarrow } \zeta$ between the sections of two given vector bundles $F_0$ and $F_1$ of respective fiber dimension $m$ and $p$. Multiplying the equations ${\cal{D}}_1 \eta = \zeta$ by $p$ test functions $\lambda$ considered as a section of the adjoint vector bundle $ad(F_1)= {\wedge}^nT^*\otimes F^*_1$ and integrating by parts, we may introduce the adjoint vector bundle $ad(F_0)={\wedge}^nT^*\otimes F^*_0$ with sections $\mu$ in order to obtain the adjoint operator $ \mu \stackrel{ad({\cal{D}}_1)}{\longleftarrow} \lambda $, writing on purpose the arrow backwards, that is from right to left. As any operator is the adjoint of another operator because $ad(ad({\cal{D}}))=\cal{D}$, we may {\it decide} to denote by $\xi \stackrel{ad({\cal{D}})}{\longleftarrow}\mu $ the generating CC of $ad({\cal{D}}_1)$ by introducing a vector bundle $E$ with sections $\xi$ and its adjoint $ad(E)={\wedge}^nT^*\otimes E^*$ with sections $\nu$. We have thus obtained the formally exact differential sequence:  \\
\[    \nu  \stackrel{ad({\cal{D}})}{\longleftarrow}  \mu \stackrel{ad({\cal{D}}_1)}{\longleftarrow} \lambda \]
and its formaly exac adjoint:    \\
\[    \xi   \stackrel{\cal{D}}{\longrightarrow} \eta \stackrel{{\cal{D}}_1}{\longrightarrow} \zeta  \]
providing a parametrization {\it if and only if} ${\cal{D}}_1$ generates the CC of ${\cal{D}}$, a situation that may not be satisfied but that we shall assume from now on because otherwise ${\cal{D}}_1$ {\it cannot} be parametrized according to the {\it double differential duality} test, for example in the case of the Einstein equations ([32],[50]) or the extension to the conformal group and other Lie groups of transformations ([33],[37- 39],[42-44]). Nevertheless, for the interested reader only, we provide the following key result on which this procedure is based (See [12],[28],[29] and [34] for more details):  \\

\noindent
{\bf THEOREM  2.1}: If $M$ is a differential module, we have the exact sequence of differential modules: \\
\[    0 \rightarrow t(M) \rightarrow M \stackrel{\epsilon}{\longrightarrow} hom_D(hom_D(M,D),D)  \]
where the map $\epsilon$ is defined by $\epsilon(m)(f)=f(m), \forall m\in M, f\in hom_D(M,D)$. Moreover, if $N$ is the differential module defined by $ad({\cal{D}})$, then $t(M) = ext^1_D(N,D)$.  \\

In order to pass to the differential module framework, let us introduce the free differential modules $D\xi\simeq D^l, D\eta \simeq D^m, D\zeta \simeq D^p$. We have similarly the adjoint free differential modules  $D\nu \simeq D^l, D\mu \simeq D^m, D\lambda \simeq D^p$, because $dim(ad(E))=dim(E)$ and $hom_D(D^m,D)\simeq D^m$. Of course, in actual practice, {\it the geometric meaning is totally different }because we have volume forms in the dual framework. We have thus obtained the formally exact sequence of differential modules:  \\
\[   D^p   \stackrel{{\cal{D}}_1}{\longrightarrow} D^m \stackrel{{\cal{D}}}{\longrightarrow} D^l  \] 
and the formally exact adjoint sequence:    \\
\[   D^p  \stackrel{ad({\cal{D}}_1)}{\longleftarrow}  D^m \stackrel{ad({\cal{D}})}{\longleftarrow} D^l   \]
The procedure with $4$ steps is as follows in the operator language:  \\

\noindent
$\bullet$  STEP $1$: Start with the formally exact {\it parametrizing sequence} already constructed by differential biduality.  We have thus $im({\cal{D}})=ker({\cal{D}}_1)$ and the corresponding differential module $M_1$ defined by ${\cal{D}}_1$ is torsion-free by assumption.  \\
$\bullet$  STEP $2$: Construct the adjoint sequence which is also formally exact by assumption.  \\
$\bullet$  STEP $3$: Find a {\it maximum} set of differentially independent CC $ad({\cal{D}}'):\mu \rightarrow {\nu}' $ among the generating CC $ad({\cal{D}} ):\mu \rightarrow  \nu$ of $ad({\cal{D}}_1)$ in such a way that $im(ad({\cal{D}}'))$ is a maximum free differential submodule of $im(ad({\cal{D}}))$ that is any element in $im(ad({\cal{D}}))$ is differentially algebraic over $im(ad({\cal{D}}'))$. \\
$\bullet$  STEP $4$:  Using differential duality, construct ${\cal{D}}'=ad(ad({\cal{D}}'))$.            \\

It remains to prove that $ {\cal{D}}_1$ generates the CC of $ {\cal{D}}' $ in the following diagram:   \\

  \[  \begin{array}{rcccccccl}
  \fbox{4} \hspace{1cm} &  &  & {\xi}' & & & &  &          \\
   &  &  &\uparrow  &\stackrel{{\cal{D}}'}{\searrow} & & &  &  \\
   &  &  & \xi  & \stackrel{{\cal{D}}}{\longrightarrow} &  \eta & \stackrel{{\cal{D}}_1}{\longrightarrow} & \zeta         &\hspace{1cm}   \fbox{1}  \\
   &  &  &  &  &  &  &  &  \\
   &  &  &  &  &  &  &  &  \\
  &  &  & \nu & \stackrel{ad({\cal{D}})}{\longleftarrow} & \mu & \stackrel{ad({\cal{D}}_1)}{\longleftarrow} & \lambda &\hspace{1cm} \fbox{2}    \\
    &  &  &  \uparrow  &  \stackrel{ad({\cal{D}}')}{\swarrow} &  &  &   \\
 \fbox{3} \hspace{1cm}   &  &  &  {\nu}'   &  &  &  &  & \\
    &  &  \swarrow & \uparrow &  &  &  &  &   \\
    &  0 &   & 0 &  &  &  &  &
  \end{array}    \] 
\vspace{2mm}  \\

\noindent
{\bf PROPOSITION  2.2}: ${\cal{D}}'$ is a minimum parametrization of ${\cal{D}}_1$.  \\

\noindent
{\it Proof}: Let us denote the number of potentials $\xi$ by $l$ (respectively ${\xi}'$ by $l'$), the number of unknowns $\eta$ by $m$ and the number of given equations $\zeta$ by $p$. As $ad({\cal{D}}')$ has no CC by construction, then $ad({\cal{D}}'): \mu \rightarrow {\nu}'$ is a formally surjective operator. On the differential module level, we have the injective operator $ ad({\cal{D}}'):D^{l'} \rightarrow D^m $ because there are no CC. Applying $hom_D(\bullet,D)$ or duality, we get an operator $D^m \rightarrow D^{l'}$ with a cokernel which is a torsion module because it has rank $l' - rk_D({\cal{D}}')= l' - rk_D(ad({\cal{D}}'))= l' - l' = 0$. \\

However, in actual practice as will be seen in the contact case, things are not so simple and we shall use the following commutative and exact diagram of differential modules based on a long $ker/coker$ long exact sequence (Compare to [35], and [46]):   \\
\[  \begin{array}{rcccccl}
0 \rightarrow ker(ad({\cal{D}})) \rightarrow & D^l &  &  \stackrel{ad({\cal{D}})}{\longrightarrow} &  & D^m &\rightarrow coker(ad({\cal{D}})) \rightarrow 0  \\  
   &   &  \searrow & &   \nearrow &  &     \\
   &   &                &  L  &   &    &   \\
    &   &   \nearrow  & \uparrow  &   \searrow   &   &   \\
    &  0  &     &  D^{l'} &   &   0   &  \\
    &   &   &     \uparrow   &  &   &  \\
    &   &   &        0       &   &   & 
\end{array}    \]
Setting $L =D^l / ker(ad({\cal{D}}))$ and introducing the biggest free differential module $D^{l'}\subseteq L$ we have $rk_D(D^{l'}) = rk_D(L) \leq rk_D(D^l) \Rightarrow  l' \leq l$, we may {\it define } the injective ({\it care}) operator $ad({\cal{D}}')$ by the composition of monomorphisms $D^{l'} \rightarrow L  \rightarrow D^m$ where the second is obtained by picking a basis of $D^{l'}$, lifting it to $D^l$ and pushing it to $D^m$ by applying $ad({\cal{D}})$. We notice that $L$ can be viewed as the differential module defined by the generating CC of $ad({\cal{D}})$ that could also be used as in ([35]).  \\

Then we have $ ad({\cal{D}}')\circ ad({\cal{D}}_1)= ad ({\cal{D}}_1 \circ {\cal{D}}')=0   \Rightarrow   {\cal{D}}_1\circ {\cal{D}}' =0  $ and thus ${\cal{D}}_1$ is surely {\it among} the CC of ${\cal{D}}' $. Therefore, the differential sequence $ {\xi}' \stackrel{{\cal{D}}'}{\longrightarrow} \eta \stackrel{{\cal{D}}_1}{\longrightarrow} \zeta $ on the operator level or the sequence $ D^p \stackrel{{\cal{D}}_1}{\longrightarrow} D^m \stackrel{{\cal{D}}'}{\longrightarrow} D^{l'} $ on the differential module level may not be exact and we can thus apply the previous Lemma. Changing slightly the notations, we have now $B=im({\cal{D}}_1)= ker({\cal{D}}) \subseteq ker({\cal{D}}')= Z$. But we have also $rk_D(B)=m - rk_D ({\cal{D}}), \,\,\, rk(Z)=m- rk_D({\cal{D}}') \Rightarrow rk_D(H)= rk_D({\cal{D}})- rk_D({\cal{D}}') = 0 $ by construction.  \\

Taking into account the previous Lemma, we may set $coim({\cal{D}}_1)=M_1\subseteq D^l$ by assumption and consider $im({\cal{D}}')=M'_1\subseteq D^{l'}$ in order to obtain the short exact sequence of differential modules $   0 \rightarrow H \rightarrow M_1  \rightarrow M'_1 \rightarrow 0 $. As $H$ is a torsion module and the differential module $M_1$ defined by ${\cal{D}}_1$ is torsion-free by assumption, the only possibility is that $H=0$ and thus $ im({\cal{D}}_1)=ker({\cal{D}}')$, that is ${\cal{D}}'$ is a minimum parametrization of ${\cal{D}}_1$ with $l' \leq l$ potentials.  \\
\hspace*{12cm}   Q.E.D.   \\

\noindent
{\bf EXAMPLE  2.3}: {\it Contact transformations}  \\
With $m=n=3, K=\mathbb{Q}(x^1,x^2, x^3) = \mathbb{Q}(x)$, we may introduce the so-called {\it contact} $1$-form $\alpha= dx^1 - x^3 dx^2$. The system of infinitesimal Lie equations defining the infinitesimal contact transformations is obtained by eliminating the factor $\rho(x)$ in the equations ${\cal{L}}(\xi) \alpha=\rho \alpha$ where ${\cal{L}}$ is the standard Lie derivative. This system is thus only generated by ${\eta}^1$ and ${\eta}^2$ below but is not involutive and one has to introduce ${\eta}^3$ defined by the first order CC:  \\
\[  \zeta \equiv  {\partial}_3{\eta}^1 - {\partial}_2{\eta}^2 - x^3{\partial}_1{\eta}^2 + {\eta}^3=0  \]
in order to obtain the following involutive system with two equations of class $3$ and one equation of class $2$, a result leading to ${\beta}^3_1=2, {\beta}^2_1=1, {\beta}^1_1=0$:  \\
\[  \left\{  \begin{array}{lclcl}
{\eta}^3 & \equiv &{\partial}_3{\xi}^3+{\partial}_2{\xi}^2 + 2x^3{\partial}_1{\xi}^2 - {\partial}_1{\xi}^1& = & 0  \\
{\eta}^2 & \equiv & {\partial}_3{\xi}^1 - x^3 {\partial}_3{\xi}^2 & = & 0 \\
{\eta}^1 & \equiv & {\partial}_2{\xi}^1 - x^3{\partial}_2{\xi}^2 +x^3{\partial}_1{\xi}^1 - (x^3)^2{\partial}_1 {\xi}^2 - {\xi}^3 & = & 0
\end{array}  
\right.  \fbox{$\begin{array}{lll}
1 & 2 & 3     \\
1 & 2 & 3  \\
1 & 2 & \bullet  
 \end{array}  $  }    \]
The characters are thus ${\alpha}^3_1=3-2=1 < {\alpha}^2_1=3-1=2, {\alpha}^1_1=3-0=3$ with sum equal to $1+2+3=6=dim(g_1)=3\times 3 -3$. In this situation, if $M$ is the differential module defined by this system or the corresponding operator ${\cal{D}}$, we know that  $rk_D(M)={\alpha}^3_1=1=3-2=rk_D(D\xi)- rk_D({\cal{D}})$. Of course, a differential trancendence basis for ${\cal{D}}$ can be the operator ${\cal{D}}':\xi \rightarrow \{{\eta}^2,{\eta}^3\}$ but, in view of the CC, we may equally choose any couple among $\{ {\eta}^1,{\eta}^2,{\eta}^3\}$ and we obtain $rk_D({\cal{D}}')=rk_D({\cal{D}})=2$ in any case, but now ${\cal{D}}'$ is formally surjective, contrary to ${\cal{D}}$. The same result can also be obtained directly from the unique CC or the corresponding operator ${\cal{D}}_1$ defining the differential module $M_1$. Finally, we have $rk_D(M_1)=3-1=2=rk_D(D\eta)-rk_D({\cal{D}}_1)$ and we check that we have indeed $rk_D(M)+rk_D(M_1)=1+2=3=rk_D(D\xi)$.  \\

It is well known that such a system can be parametrized by the injective parametrization (See [23] and [24] for more details and the study of the general dimension $n=2p + 1$):  \\
\[ - x^3 {\partial}_3\phi + \phi={\xi}^1,\,\,\,  -{\partial}_3\phi= {\xi} ^2, \,\,\, {\partial}_2\phi + x^3{\partial}_1\phi ={\xi}^3 \,\,\, \Rightarrow  \,\,\, 
{\xi}^1 - x^3{\xi}^2=\phi  \] 
It is however not so well known and quite striking that such a parametrization can be recovered idependently by using the parametrization of the differential module defined by ${\eta}^1=0$ with potentials ${\xi}^1$ and ${\xi}^2$ while setting:\\
\[  ({\xi}^1,{\xi}^2) \longrightarrow {\xi}^3= {\partial}_2{\xi}^1 - x^3{\partial}_2{\xi}^2 +x^3{\partial}_1{\xi}^1 - (x^3)^2{\partial}_1 {\xi}^2    \]
Taking into account the differential constraint ${\eta}^2  \equiv  {\partial}_3{\xi}^1 - x^3 {\partial}_3{\xi}^2  =  0 $, that is ${\xi}^2 = -  {\partial}_3({\xi}^1 - x^3 {\xi}^2)$ and substituting in ${\eta}^3=0$, we get no  additional constraint. We finally only need to modify the potentials while "{\it defining} " now 
$\phi = {\xi}^1 - x^3 {\xi}^2={\bar{\xi}}^1$ as before.  \\
The associated differential sequence is:   \\
\[  0 \rightarrow \phi \stackrel{{\cal{D}}_{-1}}{\longrightarrow} \xi \stackrel{{\cal{D}}}{\longrightarrow} \eta \stackrel{{\cal{D}}_1}{\longrightarrow} \zeta \rightarrow 0  \]
\[  0 \rightarrow 1 \longrightarrow 3  \longrightarrow 3  \longrightarrow 1  \rightarrow 0  \]
with Euler-Poincar\'{e} characteristic $1 - 3 + 3 - 1 =0$ but is {\it not} a Janet sequence because ${\cal{D}}_{-1}$ is not involutive, its completion to involution being the trivially involutive operator $j_1:\phi \rightarrow j_1(\phi)$. \\
Introducing the ring $D=K[d_1,d_2,d_3]=K[d]$ of linear differential operators with coefficients in the differential field $K$, the corresponding differential module $M\simeq D$ is projective and even free, thus torsion-free or $0$-pure, being defined by the split exact sequence of free differential modules:  \\
\[      0  \rightarrow D  \stackrel{{\cal{D}}_1}{\longrightarrow} D^3 \stackrel{{\cal{D}}}{\longrightarrow} D^3 \stackrel{{\cal{D}}_{-1}}{ \longrightarrow} D \rightarrow 0  \]
We let the reader prove as an exercise that the adjoint sequence:   \\
\[    0 \leftarrow \theta \stackrel{ad({\cal{D}}_{-1})}{\longleftarrow} \nu \stackrel{ad({\cal{D}})}{\longleftarrow} \mu \stackrel{ad({\cal{D}}_1)}{\longleftarrow} \lambda \leftarrow 0  \]
\[    0  \, \leftarrow 1  \,\, \, \, \longleftarrow  \,\,\, 3   \,\,  \longleftarrow  \,\, 3   \, \,\, \longleftarrow  \,\, \,   1   \,  \leftarrow 0  \]
starting from the Lagrange multiplier $\lambda$ is  also a split exact sequence of free differential modules.  \\

We finaly prove that the situation met for the contact structure is {\it exactly} the same as the one that we shall meet in the metric structure, namely that one can identify ${\cal{D}}_{-1}$ not with ${\cal{D}}_1$ of course but with $ad({\cal{D}}_1)$. For this, let us modify the "basis" linearly by setting $({\bar{\xi}}^1 = {\xi}^1 - x^3{\xi}^2, {\bar{\xi}}^2= {\xi}^2, {\bar{\xi}}^3={\xi}^3)$ and suppressing the bar for simplicity, we obtain the new injective parametrization:   \\
\[   \phi= {\xi}^1, \,\,\,  - {\partial}_3 \phi = {\xi}^2, \,\,\,  {\partial}_2\phi + x^3 {\partial}_1 \phi = {\xi}^3  \]
and may eliminate $\phi$ in order to consider the new involutive system, renumbering the equations through a cyclic permutation of $(1, 2, 3)$:  \\
\[  \left\{  \begin{array}{lclcl}
{\eta}^1 & \equiv &{\partial}_3{\xi}^3+{\partial}_2{\xi}^2 + x^3{\partial}_1{\xi}^2 - {\partial}_1{\xi}^1& = & 0  \\
 {\eta}^3 & \equiv & {\partial}_3{\xi}^1 + {\xi}^2 & = & 0 \\ 
 {\eta}^2 & \equiv & {\partial}_2{\xi}^1 + x^3{\partial}_1{\xi}^1 - {\xi}^3 & = & 0
\end{array}  \right.  
\fbox{$\begin{array}{lll}
1 & 2 & 3     \\
1 & 2 & 3  \\
1 & 2 & \bullet  
\end{array}  $  }    \]
with the unique first order CC defining ${\cal{D}}_1$:   \\
\[    \zeta \equiv   {\partial}_3{\eta}^2 - {\partial}_2{\eta}^3 - x^3 {\partial}_1{\eta}^3 + {\eta}^1=0   \]
Multiplying by $\lambda$ and integrating by parts, we obtain for $ad({\cal{D}}_1)$:   \\
\[   {\eta}^1 \rightarrow \lambda = {\mu}^1, \,\,\, {\eta}^2 \rightarrow  - {\partial}_3\lambda = {\mu}^2, \,\,\, {\eta}^3 \rightarrow {\partial}_2 \lambda  + x^3 {\partial}_1\lambda = {\mu}^3  \]
obtaining therefore ${\cal{D}}_{-1}= ad({\cal{D}}_1) \Leftrightarrow {\cal{D}}_1= ad({\cal{D}}_{-1})$ {\it exactly}. \\
As for ${\cal{D}}\xi=\eta$, we obtain the formal operator matrix: \\
\[   \left (  \begin{array}{ccc}
- d_1 & d_2 + x^3 d_1 & d_3  \\
d_2 + x^3d_1 & 0 & -1  \\
d_3 & 1  & 0
\end{array}  \right )  \] 
Similarly, for $ad({\cal{D}})$ we obtain the formal operator matrix:  \\
\[   \left (  \begin{array}{ccc}
 d_1 & -(d_2 + x^3 d_1) & -d_3  \\
-(d_2 + x^3d_1) & 0 & 1  \\
-d_3 & -1  & 0
\end{array}  \right )  \] 
and finally discover that $ad({\cal{D}})= - {\cal{D}}$, a striking result showing that both operators have the same CC and parametrization even though ${\cal{D}}$ is {\it not} self-adjoint.   \\    \\

\noindent
{\bf 3) EINSTEIN EQUATIONS}  \\

Linearizing the {\it Ricci} tensor ${\rho}_{ij}$ over the Minkowski metric $\omega$, we obtain the usual second order homogeneous {\it Ricci} operator $\Omega \rightarrow R$ with $4$ terms:  \\
\[  2 R_{ij}= {\omega}^{rs}(d_{rs}{\Omega}_{ij}+ d_{ij}{\Omega}_{rs} - d_{ri}{\Omega}_{sj} - d_{sj}{\Omega}_{ri})= 2R_{ji}  \]
\[ tr(R)= {\omega}^{ij}R_{ij}={\omega}^{ij}d_{ij}tr(\Omega)-{\omega}^{ru}{\omega}^{sv}d_{rs}{\Omega}_{uv}  \]
We may define the $Einstein$ operator by setting $E_{ij}=R_{ij} - \frac{1}{2} {\omega}_{ij}tr(R)$ and obtain the $6$ terms ([7]):  \\
\[ 2E_{ij}= {\omega}^{rs}(d_{rs}{\Omega}_{ij} + d_{ij}{\Omega}_{rs} - d_{ri}{\Omega}_{sj} - d_{sj}{\Omega}_{ri})
- {\omega}_{ij}({\omega}^{rs}{\omega}^{uv}d_{rs}{\Omega}_{uv} - {\omega}^{ru}{\omega}^{sv}d_{rs}{\Omega}_{uv})  \]
We have the (locally exact) differential sequence of operators acting on sections of vector bundles where the order of an operator is written under its arrow.:  \\
\[    T \underset 1{\stackrel{Killing}{\longrightarrow}} S_2T^* \underset 2{\stackrel{Riemann}{\longrightarrow}} F_1 \underset 1{\stackrel{Bianchi}{\longrightarrow}} F_2  \]
\[    n \stackrel{{\cal{D}}}{ \longrightarrow} n(n+1)/2 \stackrel{{\cal{D}}_1}{\longrightarrow} n^2(n^2-1)/12 \stackrel{{\cal{D}}_2}{\longrightarrow} n^2(n^2-1)(n-2)/24  \]
Our purpose is now to study the differential sequence onto which its right part is projecting:  \\
\[      S_2T^* \underset 2 {\stackrel{Einstein}{\longrightarrow}} S_2T^* \underset 1{\stackrel{div}{\longrightarrow}} T^*  \rightarrow 0  \]
\[  n(n+1)/2 \longrightarrow  n(n+1)/2  \longrightarrow n  \rightarrow 0  \]
and the following adjoint sequence where we have set ([13],[35],[37 - 39],[42]):   \\
 \[    Cauchy=ad(Killing),\,\,\, Beltrami=ad(Riemann), \,\,\, Lanczos=ad(Bianchi)    \]
\[   ad(T) \stackrel{Cauchy}{\longleftarrow} ad(S_2T^*)  \stackrel{Beltrami}{\longleftarrow} ad(F_1) \stackrel{Lanczos}{\longleftarrow} ad(F_2)     \]
In this sequence, if $E$ is a vector bundle over the ground manifold $X$ with dimension $n$, we may introduce the new vector bundle $ad(E)={\wedge}^nT^* \otimes E^*$ where $E^*$ is obtained from $E$ 
by inverting the transition rules exactly like $T^*$ is obtained from $T$. We have for eample $ad(T)={\wedge}^nT^*\otimes T^*\simeq {\wedge}^nT^*\otimes T \simeq {\wedge}^{n-1}T^*$ because $T^*$ is isomorphic to $T$ by using the metric $\omega$.
The $10 \times 10 $ $ Einstein$ operator matrix is induced from the $10 \times 20$ $Riemann$ operator matrix and the $10 \times 4$ $div$ operator matrix is induced from the $20 \times 20$ $Bianchi$ operator matrix. We advise the reader not familar with the formal theory of systems or operators to follow the computation in dimension $n=2$ with the $1 \times 3$ $Airy$ operator matrix, which is the formal adjoint of the $3 \times 1$ $Riemann$ operator matrix, and $n=3$ with the $6 \times 6$ $Beltrami$ operator matrix which is the formal adjoint of the $6 \times 6$ $Riemann$ operator matrix  which is easily seen to be self-adjoint up to a change of basis.\\ 
With more details, we have:  \\

\noindent
$\bullet \hspace{3mm}n=2$: The stress equations become $ d_1{\sigma}^{11}+ d_2{\sigma}^{12}=0, d_1{\sigma}^{21}+ d_2{\sigma}^{22}=0$. Their second order parametrization ${\sigma}^{11}= d_{22}\phi, {\sigma}^{12}={\sigma}^{21}= - d_{12}\phi, {\sigma}^{22}= d_{11}\phi$ has been provided by George Biddell Airy in 1863 ([2]) and is well known ([28]). We get the second order system:  \\
\[ \left\{  \begin{array}{rll}
{\sigma}^{11} & \equiv d_{22}\phi =0 \\
-{\sigma}^{12} & \equiv d_{12}\phi =0 \\
{\sigma}^{22} & \equiv d_{11}\phi=0
\end{array}
\right. \fbox{ $ \begin{array}{ll}
1 & 2   \\
1 & \bullet \\  
1 & \bullet  
\end{array} $ } \]
which is involutive with one equation of class $2$, $2$ equations of class $1$ and it is easy to check that the $2$ corresponding first order CC are just the $Cauchy$  equations. Of course, the $Airy$ function ($1$ term) has absolutely nothing to do with the perturbation of the metric ($3$ terms). With more details, when $\omega$ is the Euclidean metric, we may consider the only component:   \\
\[  \begin{array}{rcl}
tr(R)  &  = &  (d_{11} + d_{22})({\Omega}_{11} + {\Omega}_{22}) - (d_{11}{\Omega}_{11} + 2 d_{12}{\Omega}_{12}+ d_{22}{\Omega}_{22})  \\
      &  =  &  d_{22}{\Omega}_{11} + d_{11}{\Omega}_{22} - 2 d_{12}{\Omega}_{12}
      \end{array}  \]
Multiplying by the Airy function $\phi$ and integrating by parts, we discover that:   \\
      \[    Airy=ad(Riemann) \,\,\,\,  \Leftrightarrow  \,\,\,\,  Riemann = ad(Airy)   \]  \\
in the following differential sequences:
      \[ \hspace{4mm}  2 \underset 1 {\stackrel{Killing}{\longrightarrow}} 3 \underset 2 {\stackrel{Riemann}{\longrightarrow}} 1 \longrightarrow 0  \]
      \[ 0 \longleftarrow 2 \underset 1 {\stackrel{Cauchy}{\longleftarrow}} 3 \underset 2 {\stackrel{\hspace{3mm}Airy \hspace{3mm}}{\longleftarrow}} 1 \hspace{14mm}    \]
      
\noindent
$\bullet \hspace{3mm} n=3$: It is more delicate to parametrize the $3$ PD equations: \\
\[ d_1{\sigma}^{11}+ d_2{\sigma}^{12}+ d_3{\sigma}^{13}=0,\hspace{3mm} d_1{\sigma}^{21}+ d_2{\sigma}^{22}+ d_3{\sigma}^{23}=0, \hspace{3mm} d_1{\sigma}^{31}+ d_2{\sigma}^{32}+ d_3{\sigma}^{33}=0 \]
A direct computational approach has been provided by Eugenio Beltrami in 1892 ([3],[14]), James Clerk Maxwell in 1870 ([19]) and Giacinto Morera in 1892 ([14],[20]) by introducing the $6$ {\it stress functions} ${\phi}_{ij}={\phi}_{ji}$ in the {\it Beltrami parametrization}. The corresponding system:\\
\[   \left\{  \begin{array}{rll}
{\sigma}^{11} \equiv & d_{33}{\phi}_{22}+ d_{22}{\phi}_{33}-2 d_{23}{\phi}_{23}=0  \\
-{\sigma}^{12}\equiv & d_{33}{\phi}_{12}+ d_{12}{\phi}_{33}- d_{13}{\phi}_{23}- d_{23}{\phi}_{13}=0  \\
 {\sigma}^{22}\equiv & d_{33}{\phi}_{11}+ d_{11}{\phi}_{33}-2 d_{13}{\phi}_{13}=0  \\
{\sigma}^{13}\equiv & d_{23}{\phi}_{12}+ d_{12}{\phi}_{23}- d_{22}{\phi}_{13}- d_{13}{\phi}_{22} =0 \\
-{\sigma}^{23}\equiv & d_{23}{\phi}_{11}+ d_{11}{\phi}_{23}- d_{12}{\phi}_{13}- d_{13}{\phi}_{12} =0 \\
{\sigma}^{33}\equiv & d_{22}{\phi}_{11}+ d_{11}{\phi}_{22}-2 d_{12}{\phi}_{12}=0
\end{array}
\right. \fbox{ $ \begin{array}{lll}
1 & 2 & 3   \\
1 & 2 & 3  \\
1 & 2 & 3  \\
1 & 2 &  \bullet  \\
1 & 2 & \bullet  \\
1 & 2 & \bullet
\end{array} $ } \]
is involutive with $3$ equations of class $3$, $3$ equations of class $2$ and no equation of class $1$. The three characters are thus ${\alpha}^3_2=1\times 6 - 3=3< {\alpha}^2_2=2\times 6 -3=9 < {\alpha}^1_2= 3\times 6 - 0= 18$ and we have $dim(g_2)={\alpha}^1_2 + {\alpha}^2_2 + {\alpha}^3_2= 18 + 9 + 3 = 30 = dim(S_2T^*\otimes S_2T^*) - dim(S_2T^*)== 6\times 6 - 6$ ([22]). The $3$ CC are describing the stress equations which admit therefore a parametrization ... but without any geometric framework, in particular without any possibility to imagine that the above second order operator is {\it nothing else but} the {\it formal adjoint} of the {\it Riemann operator}, namely the (linearized) Riemann tensor with $n^2(n^2-1)/2=6$ independent components when $n=3$ ([35]). Breaking the canonical form of the six equations which is associated with the Janet tabular, we may rewrite the Beltrami parametrization of the Cauchy stress equations as follows, after exchanging the third row with the fourth row, keeping the ordering $\{(11)<(12)<(13)<(22)<(23)<(33)\}$:  \\
\[      \left(  \begin{array}{cccccc}
d_1& d_2 & d_3 &0 & 0 & 0 \\
 0 & d_1 &  0 & d_2 & d_3 & 0 \\
 0 & 0 & d_1 & 0 & d_2 & d_3 
\end{array}  \right)  
 \left(  \begin{array}{cccccc}
 0 & 0 & 0 & d_{33} & - 2d_{23} & d_{22} \\
 0 & - d_{33} & d_{23} & 0 & d_{13} & - d_{12}  \\
 0 & d_{23} & - d_{22} & - d_{13} & d_{12} & 0 \\
 d_{33}& 0 & - 2 d_{13} & 0 & 0 & d_{11}  \\
 - d_{23} & d_{13} & d_{12}& 0 & - d_{11} & 0 \\
 d_{22} & - 2 d_{12} & 0 & d_{11}& 0 & 0 
 \end{array} \right)  \equiv   0    \]
 as an identity where $0$ on the right denotes the zero operator. However, if  $\Omega$ is a perturbation of the metric $\omega$, the standard implicit summation used in continuum mechanics is, when $n=3$:  \\
 \[   \begin{array}{rcl}
{\sigma}^{ij}{\Omega}_{ij} & = & {\sigma}^{11}{\Omega}_{11} + 2 {\sigma}^{12}{\Omega}_{12} + 2 {\sigma}^{13}{\Omega}_{13} + {\sigma}^{22} {\Omega}_{22} + 2{\sigma}^{23}{\Omega}_{23} + {\sigma}^{33}{\Omega}_{33}  \\
   &  =  & {\Omega}_{22}d_{33}{\phi}_{11}+ {\Omega}_{33}d_{22}{\phi}_{11}- 2 {\Omega}_{23}d_{23}{\phi}_{11}+ ... \\
   &    & + {\Omega}_{23}d_{13}{\phi}_{12}+{\Omega}_{13}d_{23}{\phi}_{12}- {\Omega}_{12}d_{33}{\phi}_{12}- {\Omega}_{33}d_{12}{\phi}_{12} + ...
\end{array}  \]
because {\it the stress tensor density $\sigma$ is supposed to be symmetric}. Integrating by parts in order to construct the adjoint operator, we get:  \\
\[ \begin{array}{rcl}
 {\phi}_{11} &  \longrightarrow  &  d_{33}{\Omega}_{22} + d_{22}{\Omega}_{33} - 2 d_{23}{\Omega}_{23} \\
 {\phi}_{12} &  \longrightarrow   &  d_{13}{\Omega}_{23}+d_{23}{\Omega}_{13}-d_{33}{\Omega}_{12} - d_{12}{\Omega}_{33}
 \end{array}  \]
and so on, obtaining therefore the striking identification:  \\
\[       Riemann=ad(Beltrami)   \hspace{1cm} \Longleftrightarrow  \hspace{1cm}   Beltrami=ad(Riemann)  \]
between the (linearized ) Riemann tensor and the Beltrami parametrization. \\
Taking into account the factor $2$ involved by multiplying the second, third and fifth row by $2$, we get the new $6\times 6$ operator matrix with rank $3$:  \\
\[   \left(  \begin{array}{cccccc}
 0 & 0 & 0 & d_{33} & - 2d_{23} & d_{22} \\
 0 & - 2d_{33} & 2d_{23} & 0 & 2d_{13} & - 2d_{12}  \\
 0 & 2d_{23} & - 2d_{22} & - 2d_{13} & 2d_{12} & 0 \\
 d_{33}& 0 & - 2 d_{13} & 0 & 0 & d_{11}  \\
 - 2d_{23} & 2d_{13} & 2d_{12}& 0 & - 2d_{11} & 0 \\
 d_{22} & - 2 d_{12} & 0 & d_{11}& 0 & 0 
 \end{array} \right)     \]
clearly providing a self-adjoint operator. \\

{\it Surprisingly}, the Maxwell parametrization is obtained by keeping ${\phi}_{11}=A, {\phi}_{22}=B, {\phi}_{33}=C$ while setting ${\phi}_{12}={\phi}_{23}={\phi}_{31}=0$ in order to obtain the system:\\
\[   \left\{  \begin{array}{rl}
{\sigma}^{11} \equiv& d_{33}B + d_{22}C=0  \\
{\sigma}^{22}\equiv & d_{33}A+ d_{11}C =0 \\
- {\sigma}^{23}\equiv & d_{23}A=0  \\
{\sigma}^{33}\equiv & d_{22}A+ d_{11}B=0  \\
- {\sigma}^{13}\equiv & d_{13}B=0 \\
- {\sigma}^{12}\equiv & d_{12}C=0
\end{array}
\right. \fbox{ $ \begin{array}{lll}
1 & 2 & 3   \\
1 & 2 & 3  \\
1 & 2 & \bullet  \\
1 & 2 &  \bullet  \\
1 & \bullet & \bullet  \\
1 & \bullet & \bullet
\end{array} $ } \]
However, {\it this system may not be involutive} and no CC can be found "{\it a priori} " because the coordinate system is surely not $\delta$-regular. Indeed, effecting the linear change of coordinates ${\bar{x}}^1 = x^1, {\bar{x}}^2 = x^2, {\bar{x}}^3 = x^3 + x^2 + x^1 $ and taking out the bar for simplicity, we obtain the new involutive system:  \\
\[   \left\{  \begin{array}{l}
d_{33}C+ d_{13}C+ d_{23}C+ d_{12}C=0  \\
d_{33}B+ d_{13}B=0  \\
d_{33}A+ d_{23}A=0  \\
d_{23}C +d_{22}C - d_{13}C - d_{13}B - d_{12}C =0  \\
d_{23}A - d_{22}C + d_{13}B + 2 d_{12}C - d_{11}C=0  \\
d_{22}A + d_{22}C - 2 d_{12}C + d_{11}C + d_{11}B=0
\end{array} \right. 
\fbox{ $ \begin{array}{lll}
1 & 2 & 3   \\
1 & 2 & 3  \\
1 & 2 &  3  \\
1 & 2 &  \bullet  \\
1 &  2 & \bullet  \\
1 &  2 & \bullet
\end{array} $ } \]
and it is easy to check that the $3$ CC obtained just amount to the desired $3$ stress equations when coming back to the original system of coordinates. However, the three characters are different as we have now ${\alpha}^3_2=3 - 3=0 < {\alpha}^2_2= 2\times 3 -3=3 < {\alpha}^1_2=3 \times 3 - 0=9$ with sum equal to $dim(g_2)=6\times 3 - 6=18 -6= 12$. {\it We have thus a minimum parametrization}. \\

Again, {\it if there is a geometrical background, this change of local coordinates is hidding it totally}. Moreover, we notice that the stress functions kept in the procedure are just the ones on which ${\partial}_{33}$ is acting. The reason for such an apparently technical choice is related to very general deep arguments in the theory of differential modules that will only be explained at the end of the paper. 

The Morera parametrization is obtained similarly by keeping now ${\phi}_{23}=L, {\phi}_{13}=M, {\phi}_{12}=N$ while setting ${\phi}_{11}={\phi}_{22}={\phi}_{33}=0$, namely:   \\
\[   \left\{  \begin{array}{l}
d_{23}L=0  \\
d_{33}N - d_{13}L - d_{23}M = 0  \\
d_{13} M = 0  \\
d_{22}M - d_{23}N - d_{12}L = 0  \\
d_{11}L - d_{12}M - d_{13}N = 0  \\
d_{12} N =0
\end{array} \right.  \]
Using now the same change of coordinates as the one already done for the Maxwell parametrization, we obtain the following system with $3$ equations of (full) class $3$ and $3$ equations of class $2$ in the Pommaret basis corresponding to the Janet tabular:  \\
\[  \left\{  \begin{array}{l}
d_{33}N + d_{23}N + d_{13}N + d_{12}N = 0  \\
d_{33}M + d_{13}M = 0   \\
d_{33}L + d_{23}L = 0   \\
 (d_{23}N + d_{23}M - d_{23}L) + (d_{13}N - d_{13}M + d_{13}L) + d_{12}N = 0    \\
2 d_{23}M + ( d_{13}N - d_{13}M - d_{13}L) + d_{12}M - d_{11}L = 0 \\
d_{22}M +( d_{12}N - d_{12}M - d_{12}L) + d_{11}L = 0
\end{array}\right. 
\fbox{ $ \begin{array}{lll}
1 & 2 & 3   \\
1 & 2 & 3  \\
1 & 2 & 3  \\
1 & 2 &  \bullet  \\
1 & 2 & \bullet  \\
1 & 2 & \bullet
\end{array} $ } \]
After elementary but tedious computations ({\it that could not be avoided} !), one can prove that the $3$ CC corresponding to the $3$ dots are effectively satisfied and that they correspond to the $3$ Cauchy stress equations which are therefore parametrized. The parametrization is thus provided by an involutive operator defining a torsion module because the character ${\alpha}^3_2$ is vanishing in $\delta$-regular coordinates, just like before for the Maxwell parametrization. {\it We have thus another minimum parametrization}. Of course, such a result could not have been understood by Beltrami in $1892$ because the work of Cartan could not be adapted easily in the language of exterior forms and the work of Janet appeared only in $1920$ with no explicit reference to involution because only Janet bases are used ([11]) while the Pommaret bases have only been introduced in $1978$ ([22]). \\

On a purely computational level, we may also keeep only $\{  {\phi}_{11}, {\phi}_{12}, {\phi}_{22}\}$ and obtain the different involutive system with the same characters and, {\it in particular}, ${\alpha}^3_2=0$:  \\
\[    \left\{  \begin{array}{rll}
{\sigma}^{11} \equiv& {\partial}_{33}{\phi}_{22}=0  \\
-{\sigma}^{12}\equiv & {\partial}_{33}{\phi}_{12}=0  \\
 {\sigma}^{22}\equiv & {\partial}_{33}{\phi}_{11}=0  \\
{\sigma}^{13}\equiv & {\partial}_{23}{\phi}_{12}-{\partial}_{13}{\phi}_{22} =0 \\
-{\sigma}^{23}\equiv & {\partial}_{23}{\phi}_{11}-{\partial}_{13}{\phi}_{12} =0 \\
{\sigma}^{33}\equiv & {\partial}_{22}{\phi}_{11}+{\partial}_{11}{\phi}_{22}-2{\partial}_{12}{\phi}_{12}=0
\end{array}
\right. \fbox{ $ \begin{array}{lll}
1 & 2 & 3   \\
1 & 2 & 3  \\
1 & 2 & 3  \\
1 & 2 &  \bullet  \\
1 & 2 & \bullet  \\
1 & 2 & \bullet
\end{array} $ } \] 

So far, we have thus obtained three explicit local minimumm parametrizations of the Cauchy stress equations with $n(n-1)/2=3$ stress potentials but there may be others ([42]).  \\

\noindent 
$\bullet \hspace{3mm} n=4$: It just remains to explain the relation of the previous results with Einstein equations. The first suprising link is provided by the following technical proposition:  \\

\noindent
{\bf PROPOSITION 3.1}: The Beltrami parametrization is just described by the $Einstein$ operator when $n=3$. The same confusion existing between the $Bianchi$ operator and the $Cauchy$ operator has been made by both Einstein and Beltrami because the $Einstein$ operator and the $Beltrami$ operator are self-adjoint in arbitrary dimension $n \geq 3$, contrary to the $Ricci$ operator.     \\

\noindent
{\it Proof}: The number of components of the Riemann tensor is $dim(F_1)=n^2(n^2-1)/12$. We have the combinatorial formula $n^2(n^2-1)/12 - n(n+1)/2 = n(n+1)(n+2)(n-3)/12$ expressing that the number of components of the Riemann tensor is always greater or equal to the number of components of the Ricci tensor whenever $n>2$. Also, we have shown in many books ([22-25],[37],[38]) or papers ([42-45]) that the number of Bianchi identities is equal to $n^2(n^2-1)(n-2)/24$, that is $3$ when $n=3$ and $20$ when $n=4$. Of course, it is well known that the $div$ operator, induced as CC of the $Einstein$ operator, has $n$ components in arbitrary dimension $n\geq 3$.  \\
Accordingly, when $n=3$ we have $n^2(n^2-1)/12=n(n+1)/2=6$ and it only remains to prove that the $Einstein$ operator reduces to the $Beltrami$ operator {\it and not just to the Ricci operator}. The following formulas can be found in any textbook on general relativity:  \\

Hence the difference can only be seen when ${\omega}_{i\neq j}=0$. In our situation with $n=3$ and the Euclidean metric, we have:  \\
 \[   \begin{array}{rcl}
 2 R_{12}= 2 E_{12} &= &  (d_{11} + d_{22} + d_{33}) {\Omega}_{12}+ d_{12}({\Omega}_{11} + {\Omega}_{22} + {\Omega}_{33})  \\
   &  &- (d_{11}{\Omega}_{12} + d_{12}{\Omega}_{22} + d_{13}{\Omega}_{23})  
  - (d_{12}{\Omega}_{11}+ d_{22}{\Omega}_{12}+ d_{23}{\Omega}_{13})  \\
  & = & d_{33}{\Omega}_{12} + d_{12}{\Omega}_{33} - d_{13}{\Omega}_{23} - d_{23}{\Omega}_{13}
\end{array}  \]
\[  \begin{array}{rcl}
2 R_{11} &  = & (d_{11} + d_{22} + d_{33}){\Omega}_{11}  + d_{11} ({\Omega}_{11} + {\Omega}_{22} + {\Omega}_{33})   \\
         &   &    - 2 (d_{11}{\Omega}_{11} + d_{12}{\Omega}_{12} + d_{13}{\Omega}_{13}  \\
               &  = & (d_{22} + d_{33}) {\Omega}_{11} + d_{11}({\Omega}_{22} + {\Omega}_{33}) - 2 (d_{12}{\Omega}_{12} + d_{13} {\Omega}_{13})
\end{array}   \]
\[tr(R)= (d_{11}{\Omega}_{22} + d_{11}{\Omega}_{33} + d_{22}{\Omega}_{11} + d_{22}{\Omega}_{33} + d_{33}{\Omega}_{11} + d_{33}{\Omega}_{22}) 
    - 2 ( d_{12}{\Omega}_{12}+ d_{13}{\Omega}_{13} + d_{23}{\Omega}_{23})    \]
\[  2 E_{11}= d_{22}{\Omega}_{33} + d_{33}{\Omega}_{22} - 2 d_{23}{\Omega}_{23}    \]

In the light of modern differential geometry, comparing these results with the works of both Maxwell, Morera, Beltrami and Einstein, it becomes clear that they have been confusing the $div$ operator induced from the $Bianchi$ operator with the $Cauchy$ operator. However, it is also clear that they both obtained a possibility to parametrize the $Cauchy$ operator by means of $3$ arbitrary potential like functions in the case of Maxwell and Morera, $6$ in the case of Beltrami who explains the previous choices, and $10$ in the case of Einstein. Of course, as they were ignoring that the $Einstein$ operator was self-adjoint whenever $n\geq 3$, they did not notice that we have $Cauchy=ad(Killing)$ and they were unable to compare their reslts with the $Airy$ operator found as early as in 1870 for the same mechanical purpose when $n=2$. To speak in a rough way, the situation is similar to what could happen in the study of contact structures if one should confuse ${\cal{D}}_{-1}$ with ${\cal{D}}_1$ ([43]). Finally, using Theorem 2.1 or Proposition 2.2, we can choose a differential transcendence basis with $n(n-1)/2$ potentials that can be indexed by ${\phi}_{ij}={\phi}_{ji}$ with $ i<j$ or $1\leq i,j \leq n-1$ or even $2 \leq i,j \leq n$ when the dimension $n\geq 2$ is arbitrary (See [23] or [36] for more details on differential algebra).   \\
\hspace*{12cm}   Q.E.D.    \\

\noindent
{\bf REMARK  3.2}: In the opinion of the author of this paper who is not an historian of sciences but a specialist of mathematical physics interested by the analogy existing between {\it electromagnetism} (EM), {\it elasticity} (EL) and {\it gravitation} (GR) by using the conformal group of space-time (See [6],[24],[27],[30],[33],[40],[43-45] for related works), it is ifficult to imagine that Einstein could not have been aware of the works of Maxwell and Beltrami on the foundations of EL and tensor calculus. Indeed, not only they were quite famous when he started his research work but it must also be noticed that the Mach-Lippmann analogy ([1],[15],[16],[18]) was introduced at the same time (See [24] and [40] for more details on the field-matter couplings and the phenomenological law discovered by ... Maxwell too). The main idea is that classical variational calculus using a Lagrangian formalism must be considered as the basic scheme of a more general and powerful "{\it duality theory} " that only depends on new purely mathematical tools, namely " {\it group theory} " and " {\it differential homological algebra} " (See [25] or [38] for the theory and [42] for the applications).   \\  

The two following crucial results, still neither known nor acknowledged today, are provided by the next proposition and corresponding corollary ([36]):  \\

\noindent
{\bf PROPOSITION  3.3}: The $Cauchy$ operator can be parametrized by the formal adjoint of the $Ricci$ operator ($4$ terms) and the $Einstein$ operator ($6$ terms) is thus useless. The so-called gravitational waves equations are thus nothing else than the formal adjoint of the linearized $Ricci$ operator.  \\

{\it Proof}:  The {\it Einstein } operator $\Omega \rightarrow E$ is defined by setting $E_{ij}=R_{ij}-\frac{1}{2}{\omega}_{ij}tr(R)$ that we shall 
 write $Einstein= C \circ Ricci$ where $C:S_2T^*\rightarrow S_2T^*$ is a symmetric matrix only depending on $\omega$, which is invertible whenever $n\geq 3$. {\it Surprisingly}, we may also introduce the {\it same} linear transformation $ C:\Omega \rightarrow \bar{\Omega}=\Omega - \frac{1}{2}\omega \, tr(\Omega)$ and the unknown composite operator  ${\cal{X}}: \bar{\Omega}\rightarrow \Omega \rightarrow E$ in such a way that $Einstein= {\cal{X}}\circ C$ where ${\cal{X}}$ is defined by (See [GR], 5.1.5 p 134):  \\
 \[  2E_{ij} =   {\omega}^{rs}d_{rs}{\bar{\Omega}}_{ij} - {\omega}^{rs}d_{ri}{\bar{\Omega}}_{sj}-{\omega}^{rs}d_{sj}{\bar{\Omega}}_{ri}+{\omega}_{ij}{\omega}^{ru}{\omega}^{sv}d_{rs}{\bar{\Omega}}_{uv}  \]
Now, introducing the test functions ${\lambda}^{ij}$, we get:  \\
\[   {\lambda}^{ij} E_{ij}={\lambda}^{ij}(R_{ij} -\frac{1}{2}{\omega}_{ij} r(tr(R)=({\lambda}^{ij}-\frac{1}{2} {\lambda}^{rs}{\omega}_{rs}{\omega}^{ij})R_{ij}={\bar{\lambda}}^{ij}R_{ij}  \]
Integrating by parts while setting as usual $\Box = {\omega}^{rs}d_{rs}$, we obtain:  \\
\[   (\Box {\bar{\lambda}}^{rs}+ {\omega}^{rs}d_{ij}{\bar{\lambda}}^{ij}-{\omega}^{sj}d_{ij}{\bar{\lambda}}^{ri}- {\omega}^{ri}d_{ij}{\bar{\lambda}}^{sj}){\Omega}_{rs}  =  {\sigma}^{rs}{\Omega}_{rs}\]
Moreover, suppressing the "bar " for simplicity, we have:   \\
\[   d_r{\sigma}^{rs}={\omega}^{ij}d_{rij} {\lambda}^{rs}+{\omega}^{rs}d_{rij}{\lambda}^{ij}-
{\omega}^{sj}d_{rij}{\lambda}^{ri} - {\omega}^{ri} d_{rij}{\lambda}^{sj}=0  \]
As $Einstein$ is a self-adjoint operator (contrary to the Ricci operator), we have the identities:   \\
\[  ad(Einstein)=ad(C) \circ ad({\cal{X}}) \,  \Rightarrow  \,  Einstein = C \circ ad({\cal{X}}) \, \Rightarrow  \, ad({\cal{X}})=Ricci \, \Rightarrow \,{\cal{X}}=ad(Ricci)  \]
Indeed, $ad(C) = C$ because $C$ is a symmetric matrix and we know that $ad(Einstein)=Einstein$. Accordingly, the operator $ad(Ricci)$ parametrizes  the $Cauchy$ equations, {\it without any reference} to the $Einstein$ operator which has no mathematical origin, in the sense that it cannot be obtained by any diagram chasing.  The three terms after the {\it Dalembert} operator factorize through the divergence operator $d_i{\lambda}^{ri}$. We may thus add the {\it differential constraints} $d_i{\lambda}^{ri}=0$ {\it without any reference to a gauge transformation} in order to obtain a (minimum) {\it relative parametrization} (see [31] and [34] for details and explicit examples). When $n=4$ we finally obtain the adjoint sequences:  \\
 \[  \begin{array}{rccccl}
  &4 & \stackrel{Killing}{\longrightarrow } & 10 & \stackrel{Ricci}{\longrightarrow}& 10  \\
   &  &  &  &  &    \\
0 \leftarrow & 4 & \stackrel{Cauchy}{\longleftarrow} & 10 &\stackrel{ad(Ricci)}{\longleftarrow} & 10 
                   \end{array}  \]  
{\it without any reference} to the $Bianchi$ operator and the induced $div$ operator. \\
Finally, using Theorem 2.1 or Proposition 2.2, we may choose  a differential transcendence basis made by $\{ {\lambda}^{ij} \mid i< j \}$ or $\{ {\lambda}^{ij} \mid 1 < i,j < n-1 \}$ or even $\{ {\lambda}^{ij} \mid 2 < i,j < n \}$ when the dimension $n \geq 2$ is arbitrary (See again [23] or [36] for more details on differential algebra). \\                                                                                                                                                                                                                                                                                                                                                                                                                                                                                                                                 
\hspace*{12cm}  Q.E.D.   \\

\noindent
{\bf COROLLARY  3.4 }: The differential module $N$ defined by the $Ricci$ or the $Einstein$ operator is not torsion-free and cannot therefore be parametrized. Its 
torsion submodule is generated by the $10$ components of the linearized $Weyl$ tensor that are killed by the $Dalembert$ operator. \\

\noindent
{\it Proof}: In order to avoid using extension modules, we present the $5$ steps of the {\it double differential duality test} in this framework:  \\
Step 1: Start with the $Einstein$ operator ${\cal{D}}_1:10 \stackrel{Einstein}{\longrightarrow} 10$.  \\
Step 2: Consider its formal adjoint: $ad({\cal{D}}_1):10  \stackrel{Einstein}{\longleftarrow} 10$.  \\
Step 3: Compute the generating CC, namely the $Cauchy$ operator: $ad({\cal{D}}):4 \stackrel{Cauchy}{\longleftarrow}10$.  \\
Step 4: Consider its formal adjoint: ${\cal{D}}=ad(ad({\cal{D}})):4 \stackrel{Killing}{\longrightarrow} 10$.  \\ 
Step 5: Compute the generating CC, namely the $Riemann$ operator: ${\cal{D}}'_1:10 \stackrel{Riemann}{\longrightarrow} 20$.  \\
With a slight abuse of language, we have the direct sum $Riemann = Ricci \oplus Weyl$ with $20=10 + 10$. It follows from differential homological algebra that the $10$ additional CC in ${\cal{D}}'_1$ that are {\it not} in ${\cal{D}}_1$, are generating the torsion submodule $t(N)$ of the differentil module $N$ defined by the $Einstein$ or $Ricci$ operator. In general, if $K$ is a differential field with commuting derivations ${\partial}_1, ... ,{\partial}_n$, we way consider the ring $D=K[d_1, ... ,d_n]=K[d]$ of differential operators with coefficients in $K$ and it is know that $rk_D({\cal{D}})=rk_D(ad({\cal{D}}))$ for any operator matrix ${\cal{D}}$ with coefficients in $K$. In the present situation, as the $Minkowski$ metric has coefficients equal to $0, 1, -1$, we may choose the ground differential field to be $K= \mathbb{Q}$. Hence, there exists operators ${\cal{P}}$ and  ${\cal{Q}}$ such that we have an identity:  \\
 \[   {\cal{P}}\circ Weyl = {\cal{Q}} \circ Ricci    \]
 One may also notice that $rk_D(Einstein)=rk_D(Ricci)$ with:  \\
 \[  rk_D(Einstein)=\frac{n(n+1)}{2}- n=\frac{n(n-1) }{2}  \hspace{3mm}, \hspace{3mm}rk_D(Riemann)=\frac{n(n+1)}{2}- n=\frac{n(n-1)}{2}   \] 
The differential ranks of the Einstein and Riemann operators are thus equal, but {\it this is a pure coincidence} because $rk_D(Einstein)$ has only to do with the $div$ operator induced by contracting the $Bianchi$ operator, while $rk_D(Riemann)$ has only to do with the classical $Killing$ operator and the fact that the corresponding differential module is a torsion module because we have a Lie group of transformations having $n + \frac{n(n-1)}{2}=\frac{n(n+1)}{2}$ parameters (translations + rotations). Hence, as the $Riemann$ operator is a direct sum of the $Weyl$ operator and the $Einstein$ or $Ricci$ operator according to the previous theorem, each component of the $Weyl$ operator must be killed by a certain operator whenever the $Einstein$ or $Ricci$ equations in vacuum are satisfied. {\it It is not at all evident} that we have ${\cal{P}}=\Box$ acting on each component of the $Weyl$ operator. A direct tricky computation can be found in ([5], p 206]), ([10], exercise 7.7]) and ([37], p 95). With more details, we may start from the long exact sequence: \\
\[       0 \rightarrow \Theta \rightarrow 4  \stackrel{Killing}{\longrightarrow}10 \stackrel{Riemann}{\longrightarrow}20  \stackrel{Bianchi}{\longrightarrow} 20 \rightarrow 6 \rightarrow 0  \]
This resolution of the set of $Killing$ vector fields is {\it not} a $Janet$ sequence because the $Killing$ operator is not involutive as it is an operator of finite type with symbol of dimension $n(n-1)/2=6$ and one should need one prolongation for getting an involutive operator with vanishing second order symbol. Splitting the $Riemann $ operator we get the commutative and exact diagram:  \\
\[  \begin{array}{cccccccccl}
  &  &  &  &  0  &  &  0  &  &  0  &     \\
  &  &  & &  \downarrow & &  \downarrow  &  &  \downarrow  &  \\
  &  & 0 &  & 10 & \longrightarrow & 16 & \rightarrow & 6 &  \rightarrow 0  \\
 &  &   \downarrow &  &  \downarrow \uparrow&  & \downarrow &  & \parallel &   \\
 4 & \stackrel{Killing}{\longrightarrow} &  10  &  \stackrel{Riemann}{\longrightarrow} & 20 & \stackrel{Bianchi}{\longrightarrow} & 20 & \rightarrow & 6 & \rightarrow 0  \\
&  &  \parallel  &  &  \downarrow\uparrow &  & \downarrow &  & \downarrow &   \\
&  &  10 & \stackrel{Einstein}{\longrightarrow}&  10 &  \stackrel{div}{\longrightarrow} & 4  & \rightarrow & 0 &  \\
&  &  \downarrow &  & \downarrow & & \downarrow & & &  \\
&  &  0 & & 0 & & 0 & & & 
\end{array}   \]
Passing to the module point of view, we have the long exact sequence:  \\
\[  0 \rightarrow D^6 \longrightarrow D^{20} \stackrel{Bianchi}{\longrightarrow} D^{20} \stackrel{Riemann}{\longrightarrow} D^{10} \stackrel{Killing}{\longrightarrow}D^4 \rightarrow M \rightarrow 0  \]
which is a resolution of the $Killing$ differential module $M=coker(Killing)$ and we check that we have indeed the vanishing of the {\it Euler-Poincar\'{e} characteristic} $6 - 20 + 20 - 10 + 4 =0$.  Accordingly, we have $N'=coker(Riemann)\simeq im(Killing) \subset D^4 $ and thus $N'$ is torsion-free with $rk_D(N')=4-0=4=n$ because $rk_D(M)=0$. \\
We have the following commutative and exact diagram where $N=coker(Einstein)$:  \\
\[  \begin{array}{rcccccccccl}
  & & &  &  &  &  &  &  &  0 &  \\
  & & &  &  &  &  &  &  &  \downarrow & \\
  & &  & 0 &  & 0  & & 0 & & t(N) &    \\
  &  &  &  \downarrow &  & \downarrow  &  & \downarrow & & \downarrow &  \\
  &  0  & \longrightarrow &  D^4 &  \stackrel{div}{\longrightarrow} & D^{10} & \stackrel{Einstein}{\longrightarrow} & D^{10} & \longrightarrow & N & \rightarrow 0  \\
 & \downarrow & & \downarrow &  & \downarrow &  & \parallel & & \downarrow &  \\
    0 \rightarrow &D^6 &\longrightarrow & D^{20} & \stackrel{Bianchi}{\longrightarrow} & D^{20} & \stackrel{Riemann}{\longrightarrow} & D^{10} & \longrightarrow &  N' &  \rightarrow 0  \\
   & \parallel & & \downarrow  &  &\downarrow &  & \downarrow &  & \downarrow &   \\
   0\rightarrow & D^6 & \longrightarrow & D^{16} & \longrightarrow & D^{10}  &  & 0 &  & 0 &  \\
   & \downarrow &  & \downarrow &  & \downarrow &  &  &  &  &  \\
   &  0 &  & 0  &  &  0  &  & &  &  &
    \end{array}  \]
If $L$ is the kernel of the epimorphism $N \rightarrow N'$, it is a torsion module because $rk_D(L)=rk_D(N) - rk_D(N')=4-4=0$. We have thus $L\subseteq t(N)$ in the following commutative and exact diagram:  \\
\[ \begin{array}{rcccl}
  & 0 & & 0 & \\
    &  \downarrow &  & \downarrow &     \\
    0 \rightarrow & L & \longrightarrow & t(N) &   \\
     & \downarrow  &  & \downarrow &  \\
    0 \rightarrow & N &  = & N & \rightarrow 0  \\
    &  \downarrow &  & \downarrow &   \\
     & N' & \longrightarrow & N/t(N) & \rightarrow 0  \\
     & \downarrow & & \downarrow &  \\
     & 0  &  & 0 &   
     \end{array}  \]
where $N/t(N)$ is a torsion-free module by definition. A snake chase allows to prove that the cokernel of the monomorphism $L \rightarrow t(N)$ is isomorphic to the kernel of the induced epimorphism $N' \rightarrow N/t(N)$ and must be therefore, at the same time, a torsion module because $rk_D(L)=rk_D(t(N))=0$ and a torsion-free module because $ N'\subset D^4$, a result leading to a contradiction unless it is zero and thus $L=t(N)$. A snake chase in the previous diagram allows to exhibit the long exact connecting sequence:  \\
 \[  0 \rightarrow D^6 \longrightarrow D^{16} \longrightarrow D^{10} \longrightarrow t(N) \rightarrow 0  \]
It must be noticed that one cannot find canonical morphisms between the classical and conformal resolutions constructed similarly because we recall that, for $n=4$ ({\it only}), the CC of the Weyl operator are of order $2$ and {\it not} $1$ like the Bianchi CC for the Riemann operator (See [37] for a computer algebra checking !). However, it follows from the last theorem that the short exact sequence $0 \rightarrow D^{10} \longrightarrow D^{20} \longrightarrow D^{10} \rightarrow 0$ splits with $D^{20}\simeq D^{10} \oplus D^{10}$ but the existence of a canonical lift $D^{20} \rightarrow D^{10} \rightarrow 0$ in the above diagram does not allow to split the right column and thus $N\neq N' \oplus t(N)$ as $N'$ is not even free. Hence, one can only say that the space of solutions of Einstein equations in vacuum contains the generic solutions of the Riemann operator which are parametrized by arbitrary vector fields. As for the torsion elements, we have $t(N)=coker(D^{16} \rightarrow D^{10})$ and we may thus represent them by the components of the Weyl tensor, killed by the Dalembertian. This module interpretation may thus question the proper origin and existence of gravitational waves because the $div$ operator on the upper left part of the diagram has {\it strictly nothing to do} with the $Cauchy=ad(Killing)$ operator which cannot appear {\it anywhere} in this diagram. \\
\hspace*{12cm}  Q.E.D.  \\

\noindent
{\bf COROLLARY 3.5}: More generally, when ${\cal{D}}$ is a Lie operator of {\it finite type}, that is when $[\Theta, \Theta] \subset \Theta$ under the ordinary bracket of vector fields and $g_{q+r}=0$ for $r$ large enough, then the Spencer sequence is locally isomorphic to the tensor product of the Poincar\'{e} sequence for the exterior derivative by  a finite dimensional Lie algebra. It is thus formally exact both with its adjoint sequence. As it is known that the extension modules do not depend on the resolution used, this is the reason for which not only the Cauchy operator can be parametrized but also the {\it Cosserat couple-stress equations} $ad(D_1)$ can be parametrized by $ad(D_2)$, a result not evident at all (see [6] and [30] for explicit computations).  \\ 

\noindent
{\bf REMARK  3.6}: A similar situation is well known for the $Cauchy$-$Riemann$ equations when $n=2$. Indeed, any infinitesimal complex transformation $\xi$ must be solution of the linear first order homogeneous system ${\xi}^2_2 - {\xi}^1_1=0, {\xi}^1_2 + {\xi}^2_1=0$ of infinitesimal Lie equations though we obtain ${\xi}^1_{11} +{\xi}^1_{22}=0, {\xi}^2_{11}+{\xi}^2_{22}=0$, that is ${\xi}^1$ and ${\xi}^2$ are {\it separately} killed by the second order {\it Laplace} operator $\Delta =d_{11}+d_{22}$.  \\

\noindent
{\bf REMARK  3.7}: A similar situation is also well known for the wave equations for the EM field $F$ in electromagnetism. Indeed, starting with the first set of $Maxwell$ equations $dF=0$ and using the $Minkowski$ constitutive law in vacuum with electric constant ${\epsilon}_0$ and magnetic constant ${\mu}_0$ such that ${\epsilon}_0{\mu}_0 c^2=1$ for the seconf set of $Maxwell$ equations, a standard tricky differential elimination allows to avoid the $Lorenz$ ({\it no} "t") gauge condition for the EM potential and to obtain {\it directly} $\Box F=0$ (See [36] or [38] for the details).  \\

Using computer algebra or a direct checking with the ordering $11<12<13<22<23<33$, we obtain:
\[  E_{33}= {\omega}^{44}d_{44} {\Omega}_{33} + lower \,\,terms \]
\[  E_{23}= {\omega}^{44}d_{44} {\Omega}_{23} ...........              \]

We have therefore the following Janet tabular: \\
\noindent
\[  \fbox{ $\begin{array}{cccc}
1 & 2 & 3 & 4  \\
1 & 2 & 3 & 4  \\
1 & 2 & 3 & 4   \\
1 & 2 & 3 & 4  \\
1 & 2 &  3 & 4  \\
1 & 2 &  3 & 4 \\
1 & 2 & 3 & \bullet \\
1 & 2 & 3 & \bullet  \\
1 & 2 & 3 & \bullet  \\
1 & 2 & 3 & \bullet
\end{array} $ }   \]
we are in the position to compute the characters of the $Einstein$ operator but a similar procedure could be followed with the $Ricci$ operator. We obtain at once:  \\
\[  \begin{array}{lccccccr}
{\beta}^4_2 & = & 6 & \Rightarrow & {\alpha}^4_2 & = (10 \times 1) -6 & = 4  \\
{\beta}^3_2 & = & 4 & \Rightarrow & {\alpha}^3_2 & = (10 \times 2) -4 & = 16  \\
{\beta}^2_2 & = & 0 & \Rightarrow & {\alpha}^2_2 & = (10 \times 3) -0 & = 30  \\
{\beta}^1_2 & = & 0 & \Rightarrow & {\alpha}^1_2 & = (10 \times 4) -0 & = 40  \\
\end{array}   \]
a result leading to $dim(g_2)= {\alpha}^1_2+ {\alpha}^2_2 + {\alpha}^3_2 +{\alpha}^4_2= 90$ and 
$dim(g_3)= {\alpha}^1_2+2 {\alpha}^2_2 + 3{\alpha}^3_2 + 4{\alpha}^4_2= 164$ along with the long exact sequences:  \\
\[  0 \rightarrow g_2 \rightarrow S_2T^* \otimes S_2 T^* \rightarrow S_2 T^* \rightarrow 0  \]
\[  0 \rightarrow g_3 \rightarrow S_3T^* \otimes S_2 T^* \rightarrow T^* \otimes S_2 T^* \rightarrow T^* \rightarrow 0  \]
Now, we have by definition ${\sigma}_{\chi}(div)= \left(  \begin{array}{c}
{\chi}_1, {\chi}_2 , {\chi}_3 , {\chi}_4
\end{array}  \right)$ and:  \\
\[ {\sigma}_{\chi} \left(  \begin{array}{c}
div
\end{array}  \right)
{\sigma}_{\chi} \left(  \begin{array}{c}
Einstein
\end{array} \right)= (0, 0, 0 , 0)  \] 
As the $Einstein$ operator is self-adjoint $10\times 10$ operator matrix up to a change of basis ([32]), we obain therefore $det({\sigma}_{\chi}(Einstein))=0$ a result not evident at first sight that we shall now refine.   \\

\[  \left( \begin{array}{cccccccccc}
{\chi}_1 & {\chi}_2 & {\chi}_3 & {\chi}_4 & 0 & 0 & 0 &0 & 0 & 0 \\
0 & {\chi}_1 & 0 & 0 & {\chi}_2 & {\chi}_3 & {\chi}_4 & 0 & 0 & 0 \\
0 & 0 & {\chi}_1 & 0 & 0 & 0 & 0 & {\chi}_2 & {\chi}_3 & {\chi}_4 \\
0 & 0 & 0 & {\chi}_1 & 0 & 0 &{\chi}_2 & 0 & {\chi}_3 & {\chi}_4
\end{array} \right) \left( \begin{array}{c}
E_{11} \\
E_{12}  \\
E_{13}  \\
E_{14}  \\
E_{22}  \\
E_{23}  \\
E_{24}  \\
E_{33}  \\
E_{34}  \\
E_{44}
\end{array} \right)  \]
that must be compared with the Poincar\'{e} situation when $n=3$, namely:  \\
\[  ({\chi}_1 \, {\chi}_2 , {\chi}_3 )
  \left(   \begin{array}{ccc}
0 & - {\chi}_3 & {\chi}_2 \\
{\chi}_3 & 0 & - {\chi}_1 \\
- {\chi}_2 & {\chi}_1 & 0 
\end{array}  \right) 
= (0 , 0, 0 )   \]

\noindent
{\bf 4) CONCLUSION}  \\

 After teaching elasticity during $25$ years to high level students in some of the best french civil engineering schools, the author of this paper still keeps in mind one of the most fascinating exercises that he has set up. The purpose was to explain why a dam made with concrete is {\it always} vertical on the water-side with a slope of about $42$ degrees on the other free side in order to obtain a minimum cost and the auto-stability under cracking of the surface under water (See the introduction of [K2] for more details). Surprisingly, the main tool involved  is the approximate computation of the Airy function inside the dam. The author discovered at that time that no one of the other teachers did know that the Airy parametrization is nothing else than the adjoint of the linearized Riemann operator used as generating CC for the deformation tensor by any engineer. Being involved in General Relativity (GR) at that time, it took him $25$ years (1970-1995) to prove that the Einstein equations could not be parametrized ([26],[50]). However, nobody is a prophet in his own country and it is only now that he discovered that GR could be considered as a way to parametrize the Cauchy operator. It follows that {\it exactly the same confusion} has been done by Maxwell, Morera, Beltrami and Einstein because, in all these cases,  the operator considered is self-adjoint. As a byproduct, the variational formalism cannot allow to discover it as no engineer could have had in mind to confuse the deformation tensor with its CC in the Lagrangian used for finite elements computations. It is thus an open historical problem to know whether Einstein knew any one of the previous works done as all these researchers  were quite famous at the time he was active. In our opinion at least, the comparison of the various parametrizations described in this paper needs no comment as we have only presented {\it facts}, {\it just facts}.  \\

\noindent
{\bf 5) BIBLIOGRAPHY}  \\

\noindent
[1] Adler, F.W.: \"{U}ber die Mach-Lippmannsche Analogie zum zweiten Hauptsatz, Anna. Phys. Chemie, 22 (1907) 578-594.  \\
\noindent
[2] Airy, G.B.:  On the Strains in the Interior of Beams, Phil. Trans. Roy. Soc.London, 153 (1863) 49-80.  \\  
\noindent
[3] Beltrami, E.: Osservazioni sulla Nota Precedente, Atti della Accademia Nazionale dei Lincei Rend., 1, 5 (1892) 141-142; Collected Works, t IV . \\
\noindent
[4] Bjork, J.E. (1993) Analytic D-Modules and Applications, Kluwer (1993).  \\ 
\noindent
[5] Choquet-Bruhat, Y.: Introduction to General Relativity, Black Holes and Cosmology, Oxford University Press (2015).  \\
\noindent
[6] Cosserat, E., \& Cosserat, F.: Th\'{e}orie des Corps D\'{e}formables, Hermann, Paris, (1909).\\
\noindent
[7] Foster, J., Nightingale, J.D.: A Short Course in General relativity, Longman (1979).  \\
\noindent
[8] Gasqui, J.: Sur la R\'{e}solubilit\'{e} Locale des Equations d'Einstein, Composito Mathematica, 47, 1 (1982) 43-69.         \\
\noindent
[9] Hu, S.-T.: Introduction to Homological Algebra, Holden-Day (1968).  \\
\noindent
[10] Hughston, L.P., Tod, K.P.: An Introduction to General Relativity, London Math. Soc. Students Texts 5, Cambridge University Press 
(1990). \\
\noindent
[11] Janet, M.: Sur les Syst\`{e}mes aux D\'{e}riv\'{e}es Partielles, Journal de Math., 8 (1920) 65-151. \\
\noindent 
[12] Kashiwara, M.: Algebraic Study of Systems of Partial Differential Equations, M\'{e}moires de la Soci\'{e}t\'{e} Math\'{e}matique de France, 63 (1995) (Transl. from Japanese of his 1970 MasterÕs Thesis).  \\
\noindent
[13] Lanczos, C.: The Splitting of the Riemann Tensor, Rev. Modern Physics, 34 (1962) 379-389.  \\
\noindent
[14] Landriani, G. S.: A Note about an Exchange of Opinions between Beltrami and Morera on Elastic Equilibrium Equations, Meccanica, 52 (2017) 2801- 2806.  \\
https://doi.org/10.1007/s11012-016-0611-z   \\
\noindent
[15] Lippmann, G.: Extension du Principe de S. Carnot \`{a} la Th\'{e}orie des P\'{e}nom\`{e}nes \'{e}lectriques, C. R. Acad/ Sc. Paris, 82 (1876) 1425-1428.  \\
\noindent
[16] Lippmann, G.: \"{U}ber die Analogie zwischen Absoluter Temperatur un Elektrischem Potential, Ann. Phys. Chem., 23 (1907) 994-996. \\ 
\noindent
[17] Macaulay, F.S.: The Algebraic Theory of Modular Systems, Cambridge Tract 19, Cambridge University Press, London, 1916 (Reprinted by Stechert-Hafner Service Agency, New York, 1964).  \\
\noindent
[18]ÊMach, E.: Prinzipien der W\"{a}rmelehre, 2, Aufl., p 330, Leipzig: J.A. Barth (1900).  \\
\noindent
[19] Maxwell, J.C.: On Reciprocal Figures, Frames and Diagrams of Forces, Trans. Roy. Soc. Ediinburgh, 26 (1870) 1-40.  \\
\noindent
[20] Morera, G.: Soluzione generale della equazioni indefinite di equilibrio di un corpo continuo, Atti della Academia Nazionale dei Lincei Rend., 1, 5 (1892) 137-141 + 233 - 234.        \\
\noindent
[21] Northcott, D.G.: An Introduction to Homological Algebra, Cambridge university Press (1966).  \\
\noindent
[22] Pommaret, J.-F.: Systems of Partial Differential Equations and Lie Pseudogroups, Gordon and Breach, New York (1978); Russian translation: MIR, Moscow,(1983).\\
\noindent
[23] Pommaret, J.-F.: Differential Galois Theory, Gordon and Breach, New York (1983).\\
\noindent
[24] Pommaret, J.-F.: Lie Pseudogroups and Mechanics, Gordon and Breach, New York (1988).\\
\noindent
[25] Pommaret, J.-F.: Partial Differential Equations and Group Theory, Kluwer (1994).\\
http://dx.doi.org/10.1007/978-94-017-2539-2    \\
\noindent
[26] Pommaret, J.-F.: Dualit\'{e} Diff\'{e}rentielle et Applications, Comptes Rendus Acad\'{e}mie des Sciences Paris, S\'{e}rie I, 320 (1995) 1225-1230.  \\
\noindent
[27] Pommaret, J.-F.: Fran\c{c}ois Cosserat and the Secret of the Mathematical Theory of Elasticity, Annales des Ponts et Chauss\'ees, 82 (1997) 59-66 (Translation by D.H. Delphenich).  \\
\noindent
[28] Pommaret, J.-F.: Partial Differential Control Theory, Kluwer, Dordrecht (2001).\\
\noindent
[29] Pommaret, J.-F.: Algebraic Analysis of Control Systems Defined by Partial Differential Equations, in "Advanced Topics in Control Systems Theory", Springer, Lecture Notes in Control and Information Sciences 311 (2005) Chapter 5, pp. 155-223.\\
\noindent
[30] Pommaret, J.-F.: Parametrization of Cosserat Equations, Acta Mechanica, 215 (2010) 43-55.\\
http://dx.doi.org/10.1007/s00707-010-0292-y  \\
\noindent
[31] Pommaret, J.-F.: Spencer Operator and Applications: From Continuum Mechanics to Mathematical Physics, in "Continuum Mechanics-Progress in Fundamentals and Engineering Applications", Dr. Yong Gan (Ed.), ISBN: 978-953-51-0447--6, InTech (2012) Available from: \\
http://dx.doi.org/10.5772/35607   \\
\noindent
[32] Pommaret, J.-F.: The Mathematical Foundations of General Relativity Revisited, Journal of Modern Physics, 4 (2013) 223-239. \\
 http://dx.doi.org/10.4236/jmp.2013.48A022   \\
  \noindent
[33] Pommaret, J.-F.: The Mathematical Foundations of Gauge Theory Revisited, Journal of Modern Physics, 5 (2014) 157-170.  \\
http://dx.doi.org/10.4236/jmp.2014.55026  \\
 \noindent
[34] Pommaret, J.-F.: Relative Parametrization of Linear Multidimensional Systems, Multidim. Syst. Sign. Process., 26 (2015) 405-437.  \\
DOI 10.1007/s11045-013-0265-0   \\
\noindent
[35] Pommaret, J.-F.: Airy, Beltrami, Maxwell, Einstein and Lanczos Potentials revisited, Journal of Modern Physics, 7 (2016) 699-728. \\
\noindent
http://dx.doi.org/10.4236/jmp.2016.77068   \\
\noindent
[36] Pommaret, J.-F.: Why Gravitational Waves Cannot Exist, Journal of Modern Physics, 8 (2017) 2122-2158.  \\
https://doi.org/104236/jmp.2017.813130    \\
\noindent
[37] Pommaret, J.-F.: Deformation Theory of Algebraic and Geometric Structures, Lambert Academic Publisher (LAP), Saarbrucken, Germany (2016). A short summary can be found in "Topics in Invariant Theory ", S\'{e}minaire P. Dubreil/M.-P. Malliavin, Springer 
Lecture Notes in Mathematics, 1478 (1990) 244-254.\\
http://arxiv.org/abs/1207.1964  \\
\noindent
[38] Pommaret, J.-F.: New Mathematical Methods for Physics, Mathematical Physics Books, Nova Science Publishers, New York (2018) 150 pp.  \\
\noindent
[39] Pommaret, J.-F.: Homological Solution of the Riemann/Lanczos and Weyl/Lanczos Problems in Arbitrary Dimension, \\
https://arxiv.org/abs/1803.09610  \\
\noindent
[40]  Pommaret, J.-F.: The Mathematical Foundations of Elasticity and Electromagnetism Revisited, Journal of Modern Physics, 10 (2019) 1566-1595.     \\
 https://doi.org/10.4236/jmp.2019.1013104 (http://arxiv.org/abs/1802.02430 ) \\
\noindent
[41] Pommaret, J.-F.: Generating Compatibility Conditions and General Relativity, J. of Modern Physics, 10, 3 (2019) 371-401.  \\
\noindent
https://doi.org/10.4236/jmp.2019.103025   \\
\noindent
[42] Pommaret, J.-F.: Differential Homological Algebra and General Relativity, J. of Modern Physics, 10 (2019) 1454-1486. \\
\noindent
https://doi.org/10.4236/jmp.2019.1012097   \\
\noindent
[43] Pommaret, J.-F.: The Conformal Group Revisited, \\
\noindent
https://arxiv.org/abs/2006.03449 .  \\
\noindent
[44] Pommaret, J.-F.: Nonlinear Conformal Electromagnetism and Gravitation, \\
\noindent
https://arxiv.org/abs/2007.01710 .  \\
\noindent
[45]  Pommaret, J.-F.: A Mathematical Comparison of the Schwarzschild and Kerr Metrics, Journal of modern Physics, 11, (2020) 1672-1710.  
https://dx.doi.org/10.4236/jmp.2020.1110104     \\
\noindent
[46] Pommaret, J.-F. and Quadrat, A.: Localization and Parametrization of Linear Multidimensional Control Systems, Systems \& Control Letters, 37 (1999) 247-260.  \\
\noindent 
[47] Rotman, J.J.: An Introduction to Homological Algebra, Pure and Applied Mathematics, Academic Press (1979).  \\
\noindent
[48] Schneiders, J.-P.: An Introduction to D-Modules, Bull. Soc. Roy. Sci. Li\`{e}ge, 63, 223-295 (1994).  \\
\noindent
[49] Spencer, D.C.: Overdetermined Systems of Partial Differential Equations, Bull. Am. Math. Soc., 75 (1965) 1-114.\\
\noindent
[50] Zerz, E.: Topics in Multidimensional Linear Systems Theory, Lecture Notes in Control and Information Sciences, Springer, LNCIS 256 (2000) \\

\end{document}